\documentclass[review, authoryear]{elsarticle}

\usepackage[utf8]{inputenc}
\usepackage{hyperref}
\usepackage{cancel}
\journal{Computational Statistics and Data Analysis}

\usepackage{amsmath, amssymb, amsthm, amsfonts}
\usepackage{mathtools} 
\usepackage{mathrsfs} 
\usepackage{graphicx}
\usepackage{subcaption}
\usepackage{algpseudocode}
\usepackage[linesnumbered,ruled,vlined]{algorithm2e}
\usepackage[dvipsnames]{xcolor}
\usepackage[hang,flushmargin]{footmisc}
\usepackage{lipsum}
\usepackage{appendix}
\usepackage{xr}
\usepackage{tabularx,multirow,booktabs,setspace,caption}
\usepackage{tikz}
\usepackage[utf8]{inputenc}
\usepackage[T1]{fontenc}
\usepackage{booktabs}
\usepackage[top=2.5cm, bottom=2cm]{geometry}

\usepackage{verbatim}

\makeatletter
\def\@gobbleappendixname#1\csname thesubsection\endcsname{\Alph{section}.\arabic{subsection}}
\g@addto@macro{\appendix}{
    \renewcommand{\p@subsection}{\@gobbleappendixname}

}
\makeatother

\makeatletter
\newcommand{\algorithmfootnote}[2][\footnotesize]{%
  \let\old@algocf@finish\@algocf@finish
  \def\@algocf@finish{\old@algocf@finish
    \leavevmode\rlap{\begin{minipage}{\linewidth}
    #1#2
    \end{minipage}}%
  }%
}

\DeclareCaptionLabelSeparator*{spaced}{\\[2ex]}
\newcolumntype{Y}{>{\centering\arraybackslash}X}
  
\makeatother

\DeclareMathOperator*{\argmin}{arg\,min}
\DeclareMathOperator*{\argmax}{arg\,max}

\DeclareMathOperator{\proj}{proj}

\newtheorem{property}{Property}[section]
\newcommand{\transpose}{^\top}

\newcommand{\Bs}{\ensuremath{\mathscr{B}}\space}

\newcommand{\getc}[2]{#1_{#2}^{\phantom{}}}   
\newcommand{\getct}[2]{#1_{#2}}   
\newcommand{\getsc}[3]{#1_{#2, #3}^{\phantom{}}} 
\newcommand{\diag}[1]{\mathsf{diag}\left({#1}\right)}

\newcommand{\matrice}[1]{\mathbf{#1}}

\newcommand{\projLGLtwoOrth}[3]{\proj_{\mathcal{L_{G_{#1}}} \cap \mathcal{L}_2 \cap {#2}^\perp}\left({#3}\right)}
\newcommand{\trace}[1]{\mathsf{trace}\left({#1}\right)}

\renewcommand{\eqref}[1]{Eq.~\ref{#1}}

\newcommand{\suchthat}{\mathrm{such\;that}}
\newcommand{\st}{\mathrm{subject\;to}}

\renewcommand{\sup}{\textrm{sup}}

\newcommand{\0}{\ensuremath{\mathbf{0}}\space}
\newcommand{\1}{\ensuremath{\mathbf{1}}\space}

\newcommand{\bzero}{\ensuremath{\mathbf{0}}\space}

\newcommand{\A}{\ensuremath{\mathbf{A}}\space}
\renewcommand{\a}{\ensuremath{\mathbf{a}}\space}

\newcommand{\BDelta}{\ensuremath{\boldsymbol{\Delta}}\space}
\newcommand{\BDeltahat}{\ensuremath{\boldsymbol{\widehat\Delta}}\space}

\newcommand{\Bdeltahat}{\ensuremath{\boldsymbol{\widehat\delta}}\space}

\newcommand{\C}{\ensuremath{\mathbf{C}}\space}
\newcommand{\cc}{\ensuremath{\mathbf{c}}\space}
\newcommand{\ct}{\ensuremath{\mathbf{c}\transpose}\space}

\newcommand{\pntythes}{s}
\newcommand{\Dc}{\ensuremath{\mathbf{D_c}}\space}
\newcommand{\Dchalf}{\ensuremath{\mathbf{D_c}^{\frac{1}{2}}}\space}
\newcommand{\DcNeg}{\matrice{D}^{-1}_\mathbf{c}}
\newcommand{\DcNeghalf}{{\matrice{D}^{-\frac{1}{2}}_\mathbf{c}}}

\newcommand{\F}{\ensuremath{\mathbf{F}}\space}
\newcommand{\Ft}{\ensuremath{\mathbf{F}\transpose}\space}
\newcommand{\f}{\ensuremath{\mathbf{f}}\space}

\newcommand{\fsup}{\ensuremath{\f^\ast}\space}

\newcommand{\G}{\ensuremath{\mathbf{G}}\space}
\newcommand{\Gt}{\ensuremath{\matrice{G}\transpose}\space}
\newcommand{\g}{\ensuremath{\mathbf{g}}\space}

\newcommand{\gsup}{\ensuremath{\g^\ast}\space}

\newcommand{\grpG}{\mathcal G}
\newcommand{\grpGp}{\grpG_\p}
\newcommand{\grpGq}{\grpG_\q}
\renewcommand{\H}{\ensuremath{\mathbf{H}}\space}
\newcommand{\Ht}{\ensuremath{\mathbf{H}\transpose}\space}
\newcommand{\eye}{\ensuremath{\mathbf{I}}\space}
\renewcommand{\i}{\ensuremath{\mathbf{i}}\space}
\renewcommand{\j}{\ensuremath{\mathbf{j}}\space}

\newcommand{\Lyi}{\mathcal{L}_1}
\newcommand{\Ler}{\mathcal{L}_2}

\newcommand{\Lgrp}{\mathcal{L}_\grpG}
\newcommand{\Lone}[1]{{\left\lVert {#1} \right\rVert}_1}
\newcommand{\Ltwo}[1]{{\left\lVert {#1} \right\rVert}_2}
\newcommand{\Lgnrm}[1]{{\left\lVert {#1} \right\rVert}_{\grpG}}
\newcommand{\Blambda}{\ensuremath{\boldsymbol{\Lambda}}\space}
\newcommand{\M}{\mathbf{M}}
\newcommand{\Mhalf}{{\matrice{M}}^{\frac{1}{2}}}

\newcommand{\MNeghalf}{{\matrice{M}}^{-\frac{1}{2}}}

\renewcommand{\P}{\ensuremath{\mathbf{P}}\space}
\newcommand{\Pt}{\ensuremath{\mathbf{P}\transpose}\space}
\newcommand{\p}{\ensuremath{\mathbf{p}}\space}
\newcommand{\pt}{\ensuremath{\mathbf{p}\transpose}\space}

\newcommand{\Pjr}{\boldsymbol{\Omega}_\matrice{r}}
\newcommand{\Pjc}{\boldsymbol{\Omega}_\matrice{c}}
\newcommand{\gP}{\ensuremath{\mathbf{U}}\space} 
\newcommand{\gPt}{\ensuremath{\mathbf{U}\transpose}\space} 
\newcommand{\gp}{\ensuremath{\mathbf{u}}\space}
\newcommand{\gpt}{\ensuremath{\mathbf{u}\transpose}\space} 

\newcommand{\Q}{\mathbf{Q}}
\newcommand{\Qt}{\ensuremath{\mathbf{Q}\transpose}\space}
\newcommand{\q}{\ensuremath{\mathbf{q}}\space}
\newcommand{\qt}{\ensuremath{\mathbf{q}\transpose}\space}

\newcommand{\gQ}{\ensuremath{\mathbf{V}}\space} 
\newcommand{\gQt}{\ensuremath{\mathbf{V}\transpose}\space}
\newcommand{\gq}{\ensuremath{\mathbf{v}}\space}
\newcommand{\gqt}{\ensuremath{\mathbf{v}\transpose}\space}
\newcommand{\R}{\ensuremath{\mathbf{R}}\space}
\newcommand{\rr}{\ensuremath{\mathbf{r}}\space}

\newcommand{\Dr}{\ensuremath{\mathbf{D_r}}\space}
\newcommand{\Drhalf}{\ensuremath{\mathbf{D_r}^{\frac{1}{2}}}\space}
\newcommand{\DrNeg}{\matrice{D}^{-1}_\mathbf{r}}
\newcommand{\DrNeghalf}{{\matrice{D}^{-\frac{1}{2}}_\mathbf{r}}}
\newcommand{\DrG}{\ensuremath{\matrice{D}_{\mathbf{r},\mathcal G}}\space}

\newcommand{\DrNegG}{\matrice{D}^{-1}_{\mathbf{r},\mathcal G}}

\renewcommand{\S}{\matrice{S}}

\newcommand{\W}{\ensuremath{\mathbf{W}}\space}
\newcommand{\Whalf}{{\matrice{W}}^{\frac{1}{2}}}

\newcommand{\WNeghalf}{{\matrice{W}}^{-\frac{1}{2}}}

\newcommand{\X}{\mathbf{X}}
\newcommand{\Xt}{\mathbf{X}\transpose}
\newcommand{\x}{\mathbf{x}}

\newcommand{\Z}{\ensuremath{\mathbf{Z}}\space}
\newcommand{\Zt}{\ensuremath{\mathbf{Z}\transpose}\space}

\SetCommentSty{comm}
\SetKwComment{Comment}{$\triangleright$ }{}

\newcommand{\ha}[1]{\textcolor{black}{#1}}
\newcommand{\jcy}[1]{\textcolor{black}{#1}}
\newcommand{\vg}[1]{\textcolor{black}{#1}}
\newcommand{\ak}[1]{\textcolor{black}{#1}}

\begin{document}

\begin{frontmatter}

\title{Sparse Factor Analysis for Categorical Data with the Group-Sparse Generalized Singular Value Decomposition}

\author[address1]{Ju-Chi Yu \corref{correspondingauthors}}
\ead{Ju-Chi.Yu@camh.ca}
\author[address2]{Julie Le Borgne}
\author[address3]{Anjali Krishnan}
\author[address4]{Arnaud Gloaguen}
\author[address5,address6,address7]{Cheng-Ta Yang}
\author[address3]{Laura A. Rabin}
\author[address8]{Herv{\'e} Abdi \corref{correspondingauthors}\fnref{equalcontributors}}
\ead{herve@utdallas.edu}
\author[address9]{Vincent Guillemot \corref{correspondingauthors}\fnref{equalcontributors}}
\ead{vincent.guillemot@pasteur.fr}

\address[address1]{Campbell Family Mental Health Research Institute, Centre for Addiction and Mental Health, Toronto, Canada}
\address[address2]{Université de Lille, INSERM, CHU Lille, institut Pasteur Lille, U1167-riD-AGE, Facteurs de risque et déterminants moléculaires des maladies liées au vieillissement, Lille, France}
\address[address3]{Department of Psychology, Brooklyn College of the City University of New York, BrookBlyn, USA}
\address[address4]{Centre National de Recherche en G{\'e}nomique Humaine (CNRGH), Institut de Biologie François Jacob, CEA, Universit{\'e} Paris-Saclay, {\'E}vry, France.}
\address[address5]{Department of Psychology, National Cheng Kung University, Tainan, Taiwan}
\address[address6]{Graduate Institute Mind, Brain and Consciousness, Taipei Medical University, Taipei, Taiwan}
\address[address7]{Graduate Institute of Health and Biotechnology Law, Taipei Medical University, Taipei City, Taiwan}
\address[address8]{The University of Texas at Dallas, Richardson, TX, USA}
\address[address9]{Institut Pasteur, Universit{\'e} Paris Cit{\'e}, Bioinformatics and Biostatistics Hub, F-75015 Paris, France}

\cortext[correspondingauthors]{Corresponding authors}
\fntext[equalcontributors]{Authors contributed equally}

\begin{abstract}
Correspondence analysis, 
multiple correspondence analysis 
and
their discriminant counterparts (i.e.,
discriminant simple correspondence analysis 
and 
discriminant multiple correspondence analysis)
are methods
of choice for analyzing
multivariate 
categorical data. 
In these methods, variables are integrated 
into optimal components computed as linear combinations 
whose weights are obtained from a
generalized singular value decomposition (GSVD) 
that integrates specific 
metric constraints on the rows and columns 
of the original data matrix.
The weights of the linear combinations are, 
in turn, used to interpret the 
components, and this interpretation 
is facilitated when components are 
1) pairwise orthogonal and 
2) when the values of the weights are either 
large or small but not 
intermediate---a pattern called a simple or a sparse structure. 
To obtain such simple configurations, 
 the optimization problem solved by the GSVD is extended to include
 new constraints
 that implement component orthogonality and sparse weights.
 Because 
 multiple correspondence analysis represents qualitative
 variables by a set of binary variables,
 an additional group constraint is added to
 the optimization problem in order to sparsify the whole set
 representing one qualitative variable.
 This new algorithm---called
group-sparse GSVD (gsGSVD)---%
integrates these constraints  via an iterative projection scheme 
onto the intersection of subspaces where each subspace
implements a specific constraint.
In this paper, we expose this new algorithm and 
show how it can be adapted to the sparsification
of simple and multiple correspondence analysis, and
illustrate its applications 
with the analysis of four different data sets---each illustrating
the sparsification of a particular CA-based analysis.
\end{abstract}


\end{frontmatter}


\section{Introduction}

In contemporary research, data sets routinely
describe vast samples of observations 
(in the hundreds of thousands) with
a large number of variables (in the millions or more) that can be
quantitative, qualitative, or mixtures of these two data types.
\jcy{While principal component analysis 
(PCA) extracts components from quantitative data, 
typical component methods for qualitative data are 
correspondence analysis (CA) and 
multiple correspondence analysis (MCA, a generalization
of CA). 
Like PCA, CA and MCA extract components that summarize
the associations between 
\ha{qualitative} 
variables 
by representing each variable with a set
of binary columns corresponding to its levels
–––a coding scheme called 
\emph{disjunctive coding} or \emph{group coding}
in multivariate statistics, or \emph{one hot encoding}
in machine learning \cite[see, e.g.,][]{abdi2024}.}

Component methods create new variables 
(i.e., the \emph{components})
obtained as linear combinations of the original variables,
which optimally represent, 
on the one hand,
the similarity structure of the observations 
by maps of the component space where 
the distance between observations approximates
their similarity,
and, on the other hand, 
the correlational structure of the variables
by maps where the configuration of the variable levels 
approximates their covariance.
%
Components–––being linear 
combinations of the original variables---are often 
interpreted in terms of these variables 
and are easy to interpret when
each component is obtained from a small number of variables 
where each variable only contributes to few (ideally one) components.
Such a configuration is called 
a \textit{simple} structure 
\citep{cattell1978scientific,thurstone1935vectors}---a
concept first elaborated for factor analysis methods 
by psychometricians
\citep[specifically ``factorialists'' such as, e.g.,][]{thurstone1935vectors,cattell1978scientific}.
By contrast, when the structure is not simple, components 
are hard to interpret because the contributing 
variables
are hard to identify.
To reach a simple solution and
facilitate the interpretation of 
components,
these early factorialists 
\citep{thurstone1935vectors, thurstone1947multiple}
developed heuristic procedures such as rotation.
These procedures often simplify the interpretation
(but at the cost of losing optimality and also sometimes  orthogonality),
particularly when the data fit the factorial hypotheses
(e.g., a signal comprising a mixture of few well-defined 
orthogonal dimensions and
 independent additive Gaussian noise),
but 
large 
data sets collected without a
clear construct 
are unlikely to naturally have a simple structure
recoverable from standard rotation procedures
\citep{cattell1978scientific}.
However, rotation is rarely (if ever) used to facilitate interpretation 
in the correspondence analysis family 
as MCA originated as a descriptive method
(with data unlikely to follow the factorialist model),
but finding simple structures 
has recently become more relevant because
of the large size of modern data sets.

Another traditional approach 
\citep[see, e.g.,][]{Saporta2011,abdi2014correspondence}
selects, for further investigation, 
items (i.e., rows or columns)
whose contributions exceed the average contribution
(i.e., the inverse of the number of items),
or exceed their a priori mass when the items have different masses
(i.e., such as in CA). 
These descriptive approaches can be complemented
by some inferential procedures such as test values
\cite[a cousin of Student's $t$ statistics, 
see, e.g., ][]{lebart1984multivariate}.
More recently, computer-based resampling techniques 
(e.g., bootstrap, permutation) 
provide non-parametric equivalents---such as bootstrap ratios---to
these test values.

By contrast with these
earlier heuristics,
\emph{sparsification}---the modern approach to ``simplification''---
which originated in a multiple regression 
framework, reframes simplification
as an optimization problem 
whose goal is to
minimize the sum of squared residuals
while simultaneously minimizing 
the sum of the absolute values of the coefficients
\citep[a procedure originally called
\emph{least absolute shrinkage and selection operator}, 
better known as LASSO, see, ][]{hastie2009elements,journee2010generalized,efron2016}.
In this context, when a model is sparsified, loadings 
below a specific threshold are shrunken to zero
and so are eliminated from the model.
Often, these redundant variables provide little \emph{specific} information, 
and eliminating them from the model 
makes the prediction more reliable and 
easier to interpret \citep{efron2020prediction}.
In the early 2000's, 
LASSO based sparsification methods have been
extended to component methods,
such as principal component analysis 
(PCA, see, e.g., \citealt{jolliffe2000,jolliffe2003modified,zou2006sparse,trendafilov2014simple},
for a recent review see, \citealt{trendafilov2021}, Chapter 4).
%
Both sparsification and rotation try to find a compromise between \jcy{simplicity} and amount of
variance explained by the components.
When used for PCA, sparsification appears 
as an alternative to
rotation
but, as argued by \cite{trendafilov2015sparse},
sparsification
is to be preferred because
1) \jcy{for big data analysis, 
ease of interpretation should
be prioritized over maximizing explained variance;
and 2) concluding that small loadings obtained after a rotation 
are negligible can be misleading as }
the rotated dimensions are interpreted by imposing an arbitrary 
threshold and subjectively neglecting the loadings below this threshold 
(as illustrated by \citeauthor{cadima1995loading}, 1995).  

Recently, several LASSO-based sparsification  methods have been developed 
that be classified in
\jcy{three main approaches:} 
1) 
 the $ \Lyi $ sparsification 
\citep[i.e., the LASSO, per se,][]{tibshirani1996regression}, 
2) the $ \Lgrp $ group sparsification
\citep[i.e., the Group-LASSO,][]{yuan2006model}, 
and 
3)
 the elastic net which combines $ \Lyi $ and $ \Ler $ 
\citep{zou2005regularization}. 
These sparsification methods concurrently maximize explained variance 
while penalizing 
(i.e., reducing or eliminating) non-zero, intermediate loadings. 
In practice, sparsification is implemented by adding specific constraints 
to the component maximization problem. 

Existing sparsification methods are largely dominated by extensions 
of the penalized matrix decomposition (PMD) introduced by 
\citet{witten2009penalized}
and by sparse PCA from 
\citet{Mackey2009}. 
However, very few sparse versions of CA
and MCA have been developed. 
Notably, an early version of sparse MCA was introduced by 
\citet{Bernard2012},
whose methodology was grounded in \citet{witten2009penalized}'s PMD. 
However, 
\citet{Bernard2012}'s
approach did not rely on an explicit optimization problem 
nor did it guarantee the orthogonality of the \jcy{resulting} components. 
In a similar vein, 
\citet{Mori2016} proposed an alternative version of sparse MCA, 
based on an iterative algorithm that implements a
 decomposition of the dataset as an
optimization problem under constraints, 
including orthogonality of loadings and sparsity. 
However, this method is mostly a heuristic because
there is 
no convergence proof 
to indicate that
this algorithm solves the optimization problem
(note that as of 2023, the software associated 
with this publication is no longer available). 
In a very recent major advance,  
 \citet{Liu2023} 
presented
the first genuine sparse CA method rooted in PMD. 
In contrast with previous methods, 
 \citeauthor{Liu2023}
include a deflation step (i.e., BiOPD)
  to \ha{*** HA ***} 
minimize  
the correlation between components
(but cannot guarantee true orthogonality, because 
sparsification is not a linear problem).

These recent sparse methods 
sacrifice orthogonality for sparsity, 
because it is 
difficult 
to concurrently implement \ak{orthogonal 
components or orthogonal loadings} 
(see, \citealt{journee2010generalized,trendafilov2014simple};
exceptions are 
\citealt{trendafilov2006projected}, 
\citealt{qi2013sparse}, 
and 
\citealt{jolliffe2003modified} 
that can satisfy \ak{one of the constraints}, but not both simultaneously).
To maintain orthogonality with sparsification, 
\cite{guillemot2019constrained} 
developed a new SVD algorithm (called the constrained SVD, CSVD), 
which extends
PMD to sparsify loadings and preserves orthogonality.
Specifically, the CSVD 
imposes sparsification and orthogonality constraints by re-framing 
these constraints as convex spaces 
to solve the maximization problem 
by projecting the data onto the intersection of these spaces.

In this paper, we extend the CSVD 
to sparsify CA-related methods while maintaining the orthogonality 
of both components and factor scores. 
To do so, we extended the CSVD algorithm 
to incorporate 1) the metric matrices 
specific \jcy{to} CA and MCA and 
2) group constraints necessary for MCA
(because variables in MCA are represented by blocks of
binary variables).
Note that
these group constraints can also be applied 
\jcy{to groups of \ha{rows in CA or MCA}.}
In this paper, we show how the CSVD can be applied 
to sparsifying
CA and MCA and also to the barycentric discriminant analysis 
\citep{abdi2007discriminant}
versions of CA (DiSCA) and MCA (DiMCA).

In the following sections, we show how group sparsity 
can be used for CA-based methods, 
which properties are kept or lost, 
and how to compute optimal values 
for the sparsity parameters.
We illustrate these procedures with four examples (one per method). 
The reproducible R code of the methods 
and these examples are accessible online at 
\url{https://github.com/juchiyu/SPAFAC}.


\section{Notations}

Bold uppercase letters (e.g., $ \A $) denote a matrix, 
bold lowercase letters (e.g., $ \a $) denote a vector, 
and italic lowercase letters (e.g., $ a $)  
denote the elements of a matrix or of a vector; 
the indices for elements of a set are denoted by 
italic lowercase letters (e.g., $ n $) 
and the cardinal of a set by an italic uppercase letter (e.g., $ N $). 
The matrix $ \eye $ is the identity matrix, 
$ \1 $ is a matrix of ones, and $ \0 $ a matrix of zeros 
(the dimensions of these matrices depend upon the context 
and are assumed to be conformable). 
The transpose operation for a matrix is denoted by the superscript 
$\transpose$ (e.g., $ \A\transpose $) 
and the inverse of a square matrix (say $ \S $) 
is denoted by the superscript $ ^{-1} $ (e.g., $ \S^{-1} $). 
For a given data table, 
(except if stated otherwise)
the number of rows is denoted $ I $, 
and the number of columns  
is denoted $ J $. 
The element stored in the \textit{i}th row 
and \textit{j}th column of matrix $ \X $ is denoted $ x_{i,j} $. 
For an $I\times J$ matrix, the minimum of $ I $ and $ J $ 
is the largest possible \textit{rank}, denoted $ L $, of the matrix.
The actual rank of a matrix is denoted by $ R $ 
(with $ R \leq L $). 
The operator $\diag{}$  applied to a (square)
matrix gives a vector with all elements on the diagonal of this matrix
[i.e., $\diag{\X} = \x$, with $x_i = x_{i,i}$] and
when applied to a vector,  $\diag{}$ gives a diagonal matrix with the elements
of this vector on the diagonal  and 0s elsewhere
[i.e., $\X = \diag{\x} $, with $x_i = x_{i,i}$ and $x_{i,j} = 0$ when $i \neq j$].
The operator $ \trace{.} $ applied to a (square)
matrix gives the sum of all diagonal elements of this matrix.
The sum of squares of all elements of $ \X $ is called the 
\textit{inertia} of $ \X $ and is equal to $ \trace{\Xt\X} $ and $ \trace{\X\Xt} $. 
The Hadamar (i.e., element-wise) product 
between two matrices of same dimensions 
is denoted $\odot$ (i.e., $\mathbf{A \odot B} $).
 
The operator $ \argmax\limits_\x \left( f(\x) \right) $ 
identifies the argument $ \x $ 
that maximizes the value of $ f(\x) $. 
Similarly, the operator $ \argmin\limits_\x \left( f(\x) \right) $ 
identifies the argument $ \x $ that minimizes $ f(\x) $.
The $ \Lyi $ norm of vector $ \a $ is 
denoted by $ \Lone{\a} $, 
and  is computed as the sum of all absolute values 
of elements of  $ \a $; 
the $ \Ler $ norm of vector $ \a $, denoted by $ \Ltwo{\a} $,
is computed as $ \sqrt{\a\transpose\a} $.
In the framework of the CSVD, additional notations are needed. 
The orthogonal complement of 
$ \A $ is denoted by  the superscript 
$ ^\perp $ (e.g., $ \A^\perp $) and is the vector space 
of all the vectors orthogonal to the space spanned by $ \A $.
The space that contains all the vectors 
with an $ \Lyi $ norm inferior to $s$ is called an $ \Lyi $-ball, 
and is denoted by $ \Bs_{1}(\pntythes) $, with $ s $ 
defining the radius of the $ \Lyi $-ball.
The space that contains all the vectors
with an $ \Ler $ norm inferior to $s$ is called an $ \Ler $-ball; 
and is denoted by $ \Bs_{2}(\pntythes) $, with $ s $ 
defining the radius of the $ \Ler $-ball.
See also \ref{append:gsvd} for notation specific to the GSVD.


\section{Method}
  
CA-related methods use 
the generalized singular value decomposition (GSVD) 
to compute 
their components. 
Therefore, to sparsify these methods, 
we developed a sparsification algorithm for the GSVD,
called the sparse GSVD \citep[sGSVD;][]{Yu2023}. 
Further, as some CA-related methods analyze categorical data 
with 
qualitative variable being represented by groups
of (binary) columns,
we also extended  the sGSVD to create the group-sparse GSVD (gsGSVD)
that, in addition, performs group sparsification
where
pre-defined groups of columns or rows
are kept or eliminated together. 
In the following sections, 
we first present the optimization problem of the gsGSVD,
and the algorithm we developed to solve it. 
Next, we show how to applied the 
gsGSVD to sparsify CA, DiSCA, MCA, and DiMCA.

\subsection{The group-sparse generalized SVD (gsGSVD)}

To sparsify variables as groups of levels for CA-related methods, 
we extend the algorithm of the sGSVD to develop
the group-sparse generalized SVD
(gsGSVD), 
which sparsifies the elements of the generalized singular vectors
of a given matrix $\X$ 
 taking into account group constraints
 (on rows and columns) 
 and 
constraints imposed by a metric matrix (denoted by $\M$) for the rows  
and a metric matrix (denoted by $\W$)
for the columns (see Algorithm \ref{algo:gsGSVD}). 
Specifically,
the gsGSVD maximization problem 
is expressed as:
\begin{equation}\label{Maximization_gsGSVD_proj}
\begin{split}
&\argmax_{\getc{\p}{\ell},\getc{\q}{\ell}} 
\left( \getc{\hat\delta}{\ell} = \getct{\pt}{\ell} 
\Mhalf \X \Whalf \getc{\q}{\ell} \right)
\quad\st\quad\\[2ex]
\getc{\p}{\ell}\in
&
\begin{cases}
& \Bs_{\Ler}(1) \\
& \Bs_{\grpGp}(\getsc{\pntythes}{\p}{\ell})\\
& \getct{\P}{\ell}^\perp
\end{cases}
\quad\mathrm{,}\qquad
\getc{\q}{\ell}\in
\begin{cases}
& \Bs_{\Ler}(1) \\
& \Bs_{\grpGq}(\getsc{\pntythes}{\q}{\ell})\\
& \getct{\Q}{\ell}^\perp
\end{cases}
\end{split}
\end{equation}
with respect to three constraints, involving: 1) 
the $\Ler$-ball that normalizes the singular vectors, 
2) the orthogonal space (i.e., 
$\getct{\P}{\ell}^\perp$ 
or $\getct{\Q}{\ell}^\perp$) 
that ensures orthogonality between components, and 
3) the $\Lgrp$-ball that sparsifies the elements of the 
(generalized) singular vectors in groups.
The $\ell$th generalized pseudo-singular vector of $\X$
are then computed as
$
    \getc{\gp}{\ell} = \MNeghalf \getc{\p}{\ell} \text{ and }
    \getc{\gq}{\ell} = \WNeghalf \getc{\q}{\ell}
$.
We solve this maximization problem by modifying 
the algorithm of  the sGSVD 
(see Algorithm \ref{algo:gsGSVD}; 
cf. Algorithm \ref{algo:sGSVD} in 
\ref{append:algorithms}), 
where the three constraints are implemented with 
the Projection Onto Convex Sets (POCS) 
procedure \citep{Combettes1993pocs}. 
Specifically, the projection onto the $\Ler$-ball 
is equivalent to imposing constraints on the $\Ler$-norm, 
and the projection onto the $\Lgrp$-norm ball is equivalent
to imposing constraints on the group norm---defined as: 
$ 
\Lgnrm{\mathbf{x}} = \sum_{g = 1}^G \Ltwo{ \mathbf x_{\iota_g} } 
$ 
\citep{van2008group}. 
This group norm---%
called ``the [1,2]-group norm'' 
in \cite{van2008group}---%
is the $\Lyi$-norm of the vector containing 
the $\Ler$-norm of the sub-vectors defined by the groups.
In this paper, we call this group norm the $\Lgrp$-norm 
and its associated space the $\Lgrp$-ball, denoted $\Bs_{\grpG_{\cdot}}(\cdot)$.
\jcy{Finally}, we reorder the dimensions in a decreasing order of
the \jcy{derived} pseudo-singular values ($\getc{\hat\delta}{\ell}$). 
This \jcy{final} step is necessary because,
as is the case for the sparse SVD,
there is no guarantee that the gsGSVD 
will estimate the pseudo-singular values
in a predefined order, especially when the sparsity constraint is strong.

\begin{algorithm}[ht]
\LinesNotNumbered
 \DontPrintSemicolon
\KwData{$\X$\;
\phantom{Data:\quad}$\varepsilon$ (error), $R$ (rank, $>1$),\;
\phantom{Data:\quad}{\color{ForestGreen}$\M$ 
and $\W$ (row and column metric matrices)},\;
\phantom{Data:\quad}{\color{Maroon} 
$\getsc{\pntythes}{\p}{\ell}$ and 
$\getsc{\pntythes}{\q}{\ell}$ (sparse parameters for singular vectors)},\; 
\phantom{Data:\quad}{\color{NavyBlue} 
$\mathcal{G}_\p$, and $\mathcal{G}_\q$ 
(group arrangements of rows and columns)}}  
\KwResult{group-sparse GSVD of $\X$ }
$\widetilde{\X}$ = {\color{ForestGreen} 
$\Mhalf$}$\X${\color{ForestGreen} $\Whalf$}\; 
$\P \leftarrow \Q \leftarrow \gP \leftarrow \gQ \leftarrow \emptyset$; \Comment*[r]{$\P$ and $\Q$ store the left and right pseudo-singular vectors\\$\gP$ and $\gQ$ store the \textit{generalized} left and right pseudo-singular vectors}
\For{$\ell=1, \dots, R$}{ 
  Initialize $\p^{(0)}$ and $\q^{(0)}$
  \Comment*[r]{Initialize $\p$ and $\q$ either from SVD or randomly}
  $\hat\delta^{(0)}\leftarrow \p^{(0)}{\,}\transpose \widetilde{\X} \q^{(0)}$
    \Comment*[r]{Initialize $\hat\delta$ with $\p^{(0)}$ and $\q^{(0)}$}
  $t \leftarrow 0$\;
  \While{$\left(\Ltwo{\p^{(t+1)} - \p^{(t)}} \geq \varepsilon\right)$ \text{\rm or} $\left(\Ltwo{\q^{(t+1)} - \q^{(t)}} \geq \varepsilon\right)$
  }{
    $\p^{(t+1)} \leftarrow$ {\color{Maroon} $\proj$(} $\widetilde{\X}\q^{(t)}$, {\color{NavyBlue}
    $\Bs_{\grpGp}(\getsc{\pntythes}{\p}{\ell})\cap$}
    {\color{Maroon}$ \Bs_{\Ler}(1) \cap \P^\perp$)}\;
    $\q^{(t+1)} \leftarrow$ {\color{Maroon} $\proj$(} $\widetilde{\X}\transpose\p^{(t+1)}$, {\color{NavyBlue}
    $\Bs_{\grpGq}(\getsc{\pntythes}{\q}{\ell})\cap$}
    {\color{Maroon}$ \Bs_{\Ler}(1) \cap \Q^\perp$)}\;
    $\hat\delta^{(t+1)} \leftarrow \p^{(t+1)}{\,}\transpose \widetilde{\X} \q^{(t+1)}$\;
    $t \leftarrow t+1$ 
    \Comment*[r]{Iterate until $\p$ and $\q$ are stable}
  } 
  $\getc{\hat\delta}{\ell} 
  \leftarrow \hat\delta^{(t)}$; $\getc{\p}{\ell} \leftarrow \p^{(t)}$; 
  $\getc{\q}{\ell} \leftarrow \q^{(t)}$
  \Comment*[r]{$\getc{\p}{\ell}$ and 
  $\getc{\q}{\ell}$ 
  are the left and right pseudo-singular vectors of $\widetilde{\X}$}
  \;
  {\color{ForestGreen} 
  $\getc{\gp}{\ell}$} $\leftarrow$ 
  {\color{ForestGreen} $\MNeghalf$}$\getc{\p}{\ell}$\;
  {\color{ForestGreen} 
  $\getc{\gq}{\ell}$} $\leftarrow$ {\color{ForestGreen} $\WNeghalf$}$\getc{\q}{\ell}$
  \Comment*[r]{$\getc{\gp}{\ell}$ and 
  $\getc{\gq}{\ell}$ 
  are the left and right generalized pseudo-singular vectors of $\X$}\;
  $\P \leftarrow \left[\,\P\,\mid\,\getc{\p}{\ell}\,\right]$; 
  $\Q \leftarrow \left[\,\Q\,\mid\,\getc{\q}{\ell}\,\right]$\;
  $\gP \leftarrow \left[\,\gP\,\mid\,\getc{\gp}{\ell}\,\right]$; 
  $\gQ \leftarrow \left[\,\gQ\,\mid\,\getc{\gq}{\ell}\,\right]$\; 
}
Define $\Bdeltahat = (\hat\delta_1, \cdots, \hat\delta_R)$\;
Reorder $\Bdeltahat$ in decreasing order of its elements\;
$\BDeltahat  \leftarrow \diag{\Bdeltahat}$\;
Reorder the columns of $\gP$ and $\gQ$ according to $\Bdeltahat$\;
\caption{General algorithm of group-sparse GSVD of $\X$}
\algorithmfootnote{Note: The text colored in red is the sparsification 
constraints of the CSVD that are also used in the gsGSVD;
the text colored in green is the metric constraints 
of the GSVD that are also used in the gsGSVD; 
and the text colored in blue 
is the group constraints that are specific to the gsGSVD.}
\label{algo:gsGSVD}
\end{algorithm}

\subsection{Sparse CA and related methods}
In this paper, 
we present four sparse CA-related methods:
1) CA, 2) MCA, and 
their discriminant analysis versions 3) DiSCA and 4) DiMCA. 
These methods are all specific cases of the gsGSVD, because 
they all use metric matrices and, in some cases, 
incoporate group structures on rows or columns. 
For example, CA (and DiSCA) analyzes a contingency table and 
includes metric matrices for both rows and columns, 
whereas MCA (and DiMCA)
analyzes categorical variables
with disjunctive coding and thus, additionally, 
takes into account 
a group structure on the columns. 
In the following sections, 
we expose the algorithms for
each method as specific cases 
of the gsGSVD.

\subsubsection{CA and sparse CA (sCA): sparsification with metric matrices}
CA analyzes an $I$ by $J$ contingency table (denoted $ \A $) 
by first deriving its probability matrix 
(denoted by $ \Z $). To do so, given $ N $ observations 
(i.e., the grand total of $\A$), $\Z$ is computed as:
\begin{equation}\label{probmat_CA}
\Z = \frac{1}{\textit{N}} \,\mathbf{A}.
\end{equation}
To analyze the pattern of association between 
rows and columns, 
CA first computes the deviation from independence for $ \Z $: 
\begin{equation}\label{doubelcent_CA}
\X = \Z - \rr\ct,
\quad\mathrm{where}\quad
\quad\rr = \Z\1\quad\mathrm{and}\quad\cc = \Z\transpose\1.
\end{equation}
Here, $ \rr $, respectively $ \cc $, 
stores the sums of the rows and the columns of $ \Z $. 
Therefore, taking this deviation from independence 
is equivalent to centering both the rows and the columns of 
$ \Z $ (i.e., \textit{double-centering} $ \Z $).
CA then decomposes this double centered 
probability matrix $ \X $ into 3 matrices by the GSVD:
\begin{equation}\label{CA_GSVD}
\X = \gP \BDelta \gQt
\quad\mathrm{under\;the\;constraints}\quad
\gPt\DrNeg\gP = \gQt\DcNeg\gQ = \eye,
\end{equation}
where $\Dr$ = $\diag{\rr}$, $\Dc$ = $\diag{\cc}$.
In this GSVD, the metric 
matrix $ \DrNeg $ stores the row weights of $\X$, 
and the metric $ \DcNeg $ stores the column weights of $\X$. 
With these metric matrices, CA gives the optimal decomposition 
for a given rank of the $ \chi^2 $ 
associated to the contingency table $ \A $ 
(see \ref{append:inertia}). 
In addition, this GSVD is equivalent to the SVD of Equation 
\ref{weightX_gsvd} and solves the same maximization problem 
({cf}. 
Equation \ref{MaximizationEll_GSVD}) with 
$\M = \DrNeg$ and $\W = \DcNeg$. 
As a result, this maximization can be solved by 
\jcy{the alternating least squares (ALS) } 
Algorithm (see  Algorithm \ref{algo:gsvd} in \ref{append:algorithms}. 

In CA, each column of $\gP$ (respectively $\gQ$) 
stores the loadings for rows (respectively columns) 
to compute the row factor scores (denoted $\F$) 
and the column factor scores (denoted $\G$) of each component, where
\begin{equation}\label{f1_CA}
\F = \DrNeg\gP\BDelta = \DrNeg\X\DcNeg\gQ
\quad\mathrm{and}\quad
\G = \DcNeg\gQ\BDelta= \DcNeg\Xt\DrNeg\gP.
\end{equation}
For each component, the contributions 
of the rows and the columns are computed 
from the
squared weighted loadings:
\begin{equation}\label{ctr_CA}
    \getsc{\mathrm{\mathbf{ctr}}}{\rr}{\ell} = 
    \left( \DrNeghalf \getc{\gp}{\ell} \right) 
    \odot 
    \left( \DrNeghalf \getc{\gp}{\ell} \right)
    \quad\mathrm{and}\quad
    \getsc{\mathrm{\mathbf{ctr}}}{\cc}{\ell} = 
    \left( \DcNeghalf\getc{\gq}{\ell} \right) 
    \odot 
    \left( \DcNeghalf\getc{\gq}{\ell} \right),
\end{equation}
where the contributions of a given component sum to 1 
(i.e., $ 
\getct{\gpt}{\ell}\DrNeg\getc{\gp}{\ell} = 
\getct{\gqt}{\ell}\DcNeg\getc{\gq}{\ell} = 1 $).

For sCA, gsGSVD is used to sparsify 
the generalized singular vectors and 
gives generalized \textit{pseudo}-singular
vectors with sparsified loadings shrunken 
to 0 to facilitate the interpretation of components.
The sparsification of CA is a specific case of the gsGSVD
where the row metric matrix $\M$ equals $\DrNeg$ 
and the column metric matrix $\W$ equals $\DcNeg$. 
When the data include group structures for rows or columns
and should, therefore, be sparsified accordingly,
the group constraints can be applied 
by specifying the sparsification parameters 
$\getsc{\pntythes}{\p}{\ell}$ and 
$\getsc{\pntythes}{\q}{\ell}$ 
to include group designs of rows and of columns.
When group constraints are not included, 
the algorithm becomes the sGSVD (Algorithm 
\ref{algo:sGSVD} in \ref{append:algorithms}).

In the sCA algorithm, 
the POCS procedure projects the data onto the intersection 
of three constrained spaces:
1) the $\Ler$-ball that normalizes the singular vectors, 
2) the orthogonal space (i.e., $\getct{\P}{\ell}^\perp$ 
or $\getct{\Q}{\ell}^\perp$)
that ensures orthogonality between components, 
and
3) the $\Lyi$-ball that sparsifies the elements 
of the (generalized) singular vectors (i.e., $\p$, $\q$, $\gp$, 
and $\gq$) or the $\Lgrp$-ball if they are sparsified in groups.

From the generalized pseudo-singular vectors, 
the factor scores are computed as:
\begin{equation}\label{gFactorScores_sCA}
    \F = \DrNeg\gP\BDelta
    \quad\mathrm{and}\quad
    \G = \DcNeg\gQ\BDelta.
\end{equation}
It is worth noting that, with sparsification,
\begin{equation}\label{gFactorScores_sCA2}
    \F \neq \DrNeg\X\DcNeg\gQ
    \quad\mathrm{and}\quad
    \G \neq \DcNeg\Xt\DrNeg\gP.
\end{equation}
The loss of these properties 
is discussed in Section \ref{lostprop}.
Here, we choose to compute the factor scores
as the scaled generalized pseudo-singular vectors 
(i.e., Equation \ref{gFactorScores_sCA})
over the linear combinations of the original data 
(i.e., Equation \ref{gFactorScores_sCA2})
so that the sparsity of the generalized pseudo-singular
vectors propagates to the factor scores.

\subsubsection{Discriminant simple CA (DiSCA) and sDiSCA: 
sparsification for discriminant analysis}

The discriminant analysis of simple CA
(as opposed to \emph{multiple} CA) is used to analyze 
a contingency table with rows nested in groups 
when we aim to extract components
that best explain the inertia \textit{between}
these groups of rows. 
To do so, DiSCA performs CA on the group sums across the rows. 
Formally, DiSCA computes the matrix of group sums 
(denoted by $\A_\grpG$) with an $I \times I_\grpG $ 
indicator matrix $\H$ with 1s (and 0s otherwise) 
indicating the group belonging of the rows.
$\A_\grpG$ is therefore computed as:
\begin{equation}
\A_\grpG = \Ht\A,
\end{equation}
where the rows of $\A_\grpG$ correspond to $I_\grpG$
row groups and the columns of $\A_\grpG$ correspond 
to the $J$ levels of the categorical variable on the columns.
Similar to CA, DiSCA then computes 
the probability matrix by dividing $\A_\grpG$ 
by its total $N$:
\begin{equation}\label{probmat_DiSCA}
\Z_\grpG = \frac{1}{N} \,\A_\grpG.
\end{equation}
Next, this probability matrix $\Z_\grpG$ 
is double centered with Equation \ref{doubelcent_CA} 
to give the resulting matrix $\X_\grpG$, 
which is then decomposed by the GSVD with Equation \ref{CA_GSVD}. 
Then, we define $\DrNegG$ as the row metric matrix of $\X_\grpG$.
Finally, the row and column factors 
are computed the same way as in CA:
\begin{equation}\label{gFactorScores_DiSCA}
\F = \DrNegG\gP\BDelta = \DrNegG\X_\grpG\DcNeg\gQ
\quad\mathrm{and}\quad
\G = \DcNeg\gQ\BDelta= \DcNeg\X\transpose_{\grpG}\DrNegG\gP.
\end{equation}
The contributions of rows and columns can also be computed 
from the generalized singular vectors with Equation \ref{ctr_CA}.

\ha{In DiSCA, the row factor score $\getsc{f}{i}{\ell}$ 
represents the $i_\grpG$th group on the $\ell$th component, 
and the column factor scores $\getsc{g}{j}{\ell}$ 
represents the $j$th column on the $\ell$th component.
The individuals within each group, 
which are stored in rows of $\A$, can be projected
onto the component space as supplementary rows. }
Formally, if the $i$th row, denoted $\getct{\a\transpose}{i}$,
is projected onto the components as $\getct{\fsup}{i}$, then
\begin{equation}
    \getct{\fsup}{i} = \left(\getct{\a\transpose}{i}\1\right)^{-1}\BDelta^{-1}\G\transpose\a_i.
\end{equation}

\jcy{To sparsify DiSCA,}
sDiSCA uses the gsGSVD to sparsify the generalized singular vectors 
and computes the generalized \textit{pseudo}-singular vectors.  
Just like sCA, sDiSCA is a specific case of the gsGSVD 
where the row metric matrix $\M$ equals $\DrNegG$ 
and the column metric matrix $\W$ equals $\DcNeg$ 
with optional group constraints that can be implemented 
by specifying the sparsification parameters $\getsc{\pntythes}{\p}{\ell}$ 
and $\getsc{\pntythes}{\q}{\ell}$.
The only difference is that the row metrics are now computed
from the row groups
instead of from the individual rows.
From these generalized pseudo-singular vectors, 
the factor scores are computed in the same way as in sCA 
using Equation \ref{gFactorScores_sCA}.

Just like in sCA, 
the POCS procedure in the algorithm of sDiSCA projects 
the data onto the intersection of three constrained spaces: 
1) the $\Ler$-ball that normalizes the singular vectors, 
2) the orthogonal space 
(i.e., $\getct{\P}{\ell}^\perp$ or $\getct{\Q}{\ell}^\perp$) 
that ensures orthogonality between components, and
3) the $\Lyi$-ball that sparsifies the elements of the 
(generalized) singular vectors 
(i.e., $\p$, $\q$, $\gp$, and $\gq$) or the $\Lgrp$-ball 
if they are sparsified in groups.

\subsubsection{MCA and sparse MCA (sMCA): 
sparsification with metric matrices 
and group constraints for columns}

MCA extends CA to analyze the association between 
more than two categorical variables by using an 
\textit{indicator matrix} 
to represent each variable with a group of 0/1 columns 
\citep{greenacre1984theory,lebart1984multivariate,abdi2007multiple}. 
where each columns codes a level of this variable 
with a value of 1 indicating this level.
MCA concatenates the indicator matrices of all variables and
then uses CA to extract orthogonal components from this concatenated table. 

Formally, consider a matrix $\A$ with $ I $ observations 
and $ K $ categorical variables. 
where $ J_k $ denotes the number of levels of the
$k$th variable and with a total of 
 $ J = \sum\limits_{k = 1}^{K} J_k $ columns.
Matrix $\A$ with a sub-table structure of indicator matrices
can be expressed as
\begin{equation}\label{MCA_GroupMat}
\A = \left[\A_1 \mid \A_2 \mid \A_3 \mid \dots 
\mid \A_k \mid \dots \mid \A_K\right],
\end{equation}
where $ \A_k $ is the indicator matrix for
the $ k $th categorical variable. 
Note that 
the sum of each row of $\A$ 
equals $ K $ and the total sum of 
each $ \A_k $ equals $ I $.
Just like in CA, MCA first computes 
the probability matrix $ \Z $ by dividing $\A$ by
the sum of all it entries
(cf. Equation \ref{probmat_CA}):
\begin{equation}\label{probmat_MCA}
\Z = \frac{1}{IK} \,\A,
\end{equation}
and then double centers the probability matrix by 
(cf. Equation \ref{doubelcent_CA}):
\begin{equation}\label{doubelcent_MCA}
\X = \Z - \rr\ct,
\end{equation}
where $ \rr $ is a vector that stores the sums of each row 
(i.e., $ \rr = \Z\1 $) and $ \cc $ is a vector 
that stores the sums of each column (i.e., $ \cc = \Z\transpose\1 $). 
This double centered probability matrix $ \X $ is then analyzed by the GSVD:
\begin{equation}\label{MCA_GSVD}
\X = \gP \BDelta \gQt
\quad\st\quad
\gPt\DrNeg\gP = \gQt\DcNeg\gQ = \eye,
\end{equation}
where $\Dr = \diag{\rr}$, $\Dc = \diag{\cc}$,
and $ \BDelta $ is the diagonal matrix storing the generalized eigenvalues.
In addition, just like in CA, such a GSVD is equivalent to the SVD from 
Equation \ref{weightX_gsvd} 
and solves the same maximization problem 
({cf}. Equation \ref{MaximizationEll_GSVD})
with $\M = \DrNeg$ and $\W = \DcNeg$. 
As a result, this maximization 
can be solved by the same Algorithm \ref{algo:gsvd} as for CA. 
Finally, the row and column factors are computed as in CA:
\begin{equation}\label{gFactorScores_MCA}
\F = \DrNeg\gP\BDelta = \DrNeg\X\DcNeg\gQ
\quad\mathrm{and}\quad
\G = \DcNeg\gQ\BDelta= \DcNeg\Xt\DrNeg\gP.
\end{equation}
In MCA, the row factor scores $\getsc{\f}{i}{\ell}$ 
represents the $i$th observation on the $\ell$th component,
and the column factor scores $\getsc{\g}{j_k}{\ell}$ 
represent the $j$th level of the $k$th variable on the $\ell$th component.

To sparsify MCA, we sparsify the left and right generalized
singular vectors of MCA and derive 
\textit{generalized pseudo-singular vectors}
where the sparse loadings are shrunken to 0. 
The sparsification of MCA is a specific case of the gsGSVD 
(Algorithm \ref{algo:gsGSVD})\jcy{, similar to sCA, 
with the row metric matrix $\DrNeg$,
the column metric matrix $\DcNeg$,
and a (\textit{non-optional}) 
group constraint $\grpGq$ on the columns.}
This group constraint ensures that 
when the $k$th variable is kept or sparsified, 
all the $J_k$ columns associated 
with its levels kept or sparsified together. 
Such group sparsification 
can be optionally imposed \jcy{on the rows (i.e., the observations)} by specifying $\grpGp$
when the rows are expected to be kept or dropped according to their groups.

In the gsGSVD algorithm used for sMCA, 
the POCS procedure projects the data onto the intersection 
of three constrained spaces: 
1) the $\Ler$-ball that normalizes the singular vectors, 
2) the orthogonal space 
(i.e., $\getct{\P}{\ell}^\perp$ or $\getct{\Q}{\ell}^\perp$) 
that ensures orthogonality, and 
3) the $\Lgrp$-ball that sparsifies the elements of the right,
and sometimes also the left, 
\jcy{singular vectors (i.e., $\p$ and $\q$) 
and generalized singular vectors (i.e., $\gp$ and $\gq$)} in groups.

\subsubsection{Discriminant MCA (DiMCA) and sDiMCA: sparsification 
for discriminant analysis with group constraints for columns}

The discriminant version of MCA (DiMCA)
 analyzes
a matrix with 
variables codeed as for MCA
with rows  nested in groups
and 
extracts components that best explain 
the inertia \textit{between groups}. 
To do so, DiMCA performs MCA on the group sums 
across the rows of the concatenated table. 

Formally, consider the same $I \times J$ data set
$\A$ as used for MCA (see Equation \ref{MCA_GroupMat}),
if the observation groups are represented 
by an $ I \times I_\grpG$ indicator matrix 
$\H$, the matrix of group sums across 
the observations (denoted by $\A_\grpG$) is computed as:
\begin{equation}
\A_\grpG = \Ht\A,
\end{equation}
where the rows of $\A_\grpG$ correspond to
\jcy{the $I_\grpG$ observation groups} 
and the columns of $\A_\grpG$ correspond 
to the levels of all $K$ categorical variables 
(i.e., $J = \displaystyle\sum_{1}^{K} J_k$).
Just like in MCA, the double centered probability 
matrix \jcy{$\X_\grpG$ is obtained
from the probability matrix $\Z_\grpG$ as} 
(cf. Equations \ref{probmat_MCA} and \ref{doubelcent_MCA}):
\begin{equation}\label{probmat_DiMCA}
\Z_\grpG = \frac{1}{IK} \,\A_\grpG
\quad\mathrm{and}\quad
\X_\grpG = \Z_\grpG - \rr\ct,
\end{equation}
where $ \rr $ is a vector that stores the row sums 
(i.e., $ \rr = \Z_\grpG\1 $) and $ \cc $ 
is a vector that stores column sums 
(i.e., $ \cc = \Z_\grpG\transpose\1 $). 
$ \X_\grpG $ is then analyzed by the GSVD:
\begin{equation}\label{DiMCA_GSVD}
\X_\grpG = \gP \BDelta \gQt
\quad\st\quad
\gPt\DrNegG\gP = \gQt\DcNeg\gQ = \eye,
\end{equation}
where $\DrG = \diag{\rr}$, $\Dc = \diag{\cc}$, 
and $ \BDelta $ is the diagonal matrix 
of the generalized eigenvalues. 
In addition, just like in CA and MCA, this GSVD is equivalent 
to the SVD of Equation \ref{weightX_gsvd} 
and solves the same maximization problem 
(cf. Equation \ref{MaximizationEll_GSVD}) 
with $\M = \DrNegG$ and $\W = \DcNeg$. 
As a result, this maximization can also be solved by Algorithm \ref{algo:gsvd}. 
Finally, the row and column factors are computed 
as for CA and MCA:
\begin{equation}\label{gFactorScores_DiMCA}
\F = \DrNegG\gP\BDelta = \DrNegG\X\DcNeg\gQ
\quad\mathrm{and}\quad
\G = \DcNeg\gQ\BDelta= \DcNeg\Xt\DrNegG\gP.
\end{equation}
In DiMCA, the row factor scores $\getsc{\f}{i}{\ell}$
represent the $i_\grpG$th observation group on the $\ell$th component, 
and the column factor scores $\getsc{\g}{j_k}{\ell}$ 
represent the $j$th level of the $k$th variable on the $\ell$th component.

Finally, the original observations can be projected 
onto the component space as supplementary elements. 
Formally,  the $i$th observation denoted by $\a_i$ 
(i.e., the $i$th row of $\A$) 
is projected onto the components 
to obtain $\getct{\fsup}{i}$ computed as
\begin{equation}
    \getct{\fsup}{i} = \left(
    \getct{\a\transpose}{i}
    \1\right)^{-1}
    \getct{\a\transpose}{i}\G\BDelta^{-1}.
\end{equation}

To sparsify DiMCA, we sparsify the left and right generalized 
singular vectors of DiMCA and 
derive \textit{generalized pseudo-singular vectors}
where the sparse loadings are shrunken to 0.
Just like the other sparse methods, 
sDiMCA is a specific case of the gsGSVD (Algorithm \ref{algo:gsGSVD}) 
with row metric matrix $\DrNegG$, 
column metric matrix $\DcNeg$, 
and the group constraint $\grpGq$ (on the columns).
Here, the POCS procedure works \jcy{the same as} in sMCA.

\subsection{Sparsification: Lost and found properties}\label{lostprop}

Because of its specific preprocessing steps and metric constraints, 
CA (and therefore MCA, DiSCA, and DiMCA) 
possesses six essential properties: 
transition formulas, supplementary projections, 
asymmetric projection, distributional equivalence, 
barycentric projection, and embedded 
solutions \jcy{\citep[see][and \ref{append:properties} for details]{lebart1984multivariate}}. 
As sparsification introduces non-differentiable 
constraints into the optimization problem, 
the projection operators used for 
obtaining sparse solutions are non-linear. 
This non-linearity leads to a situation 
where the properties, 
which
predominantly depend on linearity 
in standard CA-related methods, 
are either partially retained 
with minor modifications or completely lost. 
In this section, we identify these lost properties 
and evaluate possible solutions to restore them.

\begin{property}
Transition formulas: from rows to columns and back
\end{property}
The transition formulas allow row factor scores ($\getc{\f}{\ell}$) 
to be derived from the column factor scores ($\getc{\g}{\ell}$), 
and vice versa
(see, \citealt{escofier1979representation},
the original transition formulas are also described in 
Property \ref{prop:transitionformula}). 
However, the original transition formulas no longer 
work with the gsGSVD because 
these formulas are linear projections
whereas the projecting operator in Algorithm \ref{algo:gsGSVD} 
is not linear. 
With such a non-linear projecting operator, 
we developed new transition formulas 
(with slight modifications from the original ones) 
that integrate the same non-linear projection as follows:
\begin{equation}\label{transition_sF}
    \begin{split}
        \getc{\f}{\ell} & = \DrNeg\getc{\gp}{\ell}\getc{\delta}{\ell} \\
      & = \DrNeg \cdot \Drhalf
      \projLGLtwoOrth{\gp}{\P}{\DrNeghalf\X\DcNeg
      \getc{\gq}{\ell}}\getc{\delta}{\ell} \\
      & = \DrNeghalf\projLGLtwoOrth{\gp}{\P}
      {\DrNeghalf\X\DcNeg\getc{\gq}{\ell}}
      \getc{\delta}{\ell},
    \end{split}
\end{equation}

\begin{equation}\label{transition_sG}
    \begin{split}
        \getc{\g}{\ell} & = \DcNeg\getc{\gq}{\ell}
        \getc{\delta}{\ell} \\
      & = \DcNeg\cdot\Dchalf\projLGLtwoOrth{\gq}
      {\Q}{\DcNeghalf\Xt\DrNeg
      \getc{\gp}{\ell}}\getc{\delta}{\ell} \\
      & = \DcNeghalf\projLGLtwoOrth{\gq}{\Q}
      {\DcNeghalf\Xt\DrNeg\getc{\gp}{\ell}}
      \getc{\delta}{\ell}.
    \end{split}    
\end{equation}
When $\ell = 1$, $\P^\perp = \Q^\perp = \bzero$, 
and when $\ell > 1$, 
$\P^\perp = \left[ \getc{\p}{1}\mid\dots\mid
\getc{\p}{\ell-1} \right]$ 
and $\Q^\perp = \left[ \getc{\q}{1}\mid\dots\mid
\getc{\q}{\ell-1} \right]$. 

\ha{
These transition formulas for the gsGSVD 
are valid when (1) all the constraints are satisfied,
and (2) the estimated pseudo-singular vectors
$\P$ and $\Q$ are in the order in which 
they were estimated by the ALS algorithm.
Specifically, in our algorithm, 
the number of dimensions is specified, 
and the dimensions are reordered after all 
requested ones are estimated to give 
the pseudo-singular values in descending order
\jcy{. This reordering step }is necessary 
because there is no theoretical guarantee 
that Algorithm \ref{algo:gsGSVD}
estimates pseudo-singular values in 
\jcy{a descending} order.
Remark (2) is especially important, 
because this order influences
the definition of $\P^\perp$ and $\Q^\perp$ 
for a given $\ell$.
As a result, the number of dimensions 
is one of the hyperparameters that can also 
lead to different results and should be 
evaluated and optimized (see Section \ref{section:SI} for details).
}

\begin{property}
Supplementary projections.
\end{property}

In CA/MCA/DiSCA/DiMCA, 
rows and columns can be projected 
as \textit{supplementary elements}
using  equations from Property \ref{prop:supplementaryprojection}.
Even though the equations for supplementary projection 
are directly derived from the transition formulas, 
the projection of supplementary elements no longer works
as in the regular versions of CA-related methods.
The new transition formulas are defined in Equations \ref{transition_sF} 
and \ref{transition_sG}.
When applied to a subset of existing rows or columns of $\X$,
the transition formulas are not guaranteed to recover the 
correct sparse factor scores given 
by the full sparse decomposition of $\X$;
When applied to new data, we observed that the transition formulas could give
incoherent results;
specifically
 new data numerically close to existing
data could be projected very far away from each other.
\jcy{This} problem 
occurs because the
non-linear projector 
requires the full set of  data 
(i.e., the data matrix $\X$).
\jcy{Because} $ \X $ is not fully decomposed with sparsification, 
there might not exist $ \getc{\gp}{R} $ and $ \getc{\gq}{R} $ 
for all $ R $ dimensions to perfectly reconstruct $ \X $.

To mitigate this issue, we
define a rank-$R$ approximation of $\X$ by $ \widehat{\X} $: 
\begin{equation}\label{eq:xhat}
    \widehat{\X} = \gP\widehat\BDelta\gQ\transpose,
\end{equation}
where $\gP$, $\widehat\BDelta$, $\gQ$ are given 
by Equation \ref{algo:gsGSVD}.
Based on $\widehat{\X}$, we compute two projector matrices 
for the transition formula: 
1) one projects $ \X $ onto the space of $ \widehat{\X} $ 
(denoted by $ \Pjr $) and
2) the other one projects $ \Xt $ onto the space of 
$ \widehat{\X}\transpose $ 
(denoted by $ \Pjc $). 
These linear projector matrices thus approximate 
the non-linear projectors and are obtained as:
\begin{equation}\label{proj_gsGSVD}
\begin{split}
\Pjr &= (\X\Xt)^{-1}\X(\widehat{\X}\transpose),\\
\Pjc &= (\Xt\X)^{-1}\Xt(\widehat{\X}),
\end{split}
\end{equation}
which is equivalent to using a pseudo-inverse as a projector.

Formally, the factor score of a supplementary row $\getc{\i}{\sup}$ 
(denoted $\getct{\fsup}{\sup}$) is computed 
from Equation \ref{transition_sF} as:
\begin{equation}
    \getct{\fsup}{\sup} = 
            \left(\getct{\i\transpose}{\sup}\1\right)^{-1}
            \getct{\i\transpose}{\sup}\Pjc\G\BDelta^{-1}.
\end{equation}
The factor score of a supplementary column $\getc{\j}{\sup}$ 
(denoted $\getct{\gsup}{\sup}$) is computed 
from Equation \ref{transition_sG} as:
\begin{equation}
    \getct{\gsup}{\sup} = 
    \left(\1\transpose\getc{\j}{\sup}\right)^{-1}
    \getct{\j\transpose}{\sup}\Pjr\F\BDelta^{-1}.
\end{equation}
However, when there is no linear transformation from the data in $ \X $ 
to the component space of $ \widehat{\X} $, 
these new formulas for supplementary projections will provide only an approximation.

\begin{property}
The asymmetric projection.
\end{property}

The asymmetric projections of rows or columns 
give factor scores with unitary inertia (Property \ref{prop:asymmetricprojection}). 
With sparsification, the asymmetric projections 
can be derived from the new transition formulas 
and are computed differently from Equations 
\ref{gAsymF_CA} and \ref{gAsymG_CA}. 
The asymmetric projections of the rows
$\getct{\widetilde{\F}^\ast}{\ell}$ 
are computed as:
\begin{align}\label{gAsymF_sCA}
\widetilde{\F}^\ast = \DrNeg\gP
\quad\mathrm{with\;its\;inertia\;}\quad
{\widetilde{\F}}^{\ast}{}\transpose\Dr\widetilde{\F}^{\ast} = \eye.
\end{align}
Similarly, the asymmetric projections of the columns 
$\getct{\widetilde{\G}^\ast}{\ell}$ are computed as:
\begin{align}\label{gAsymG_sCA}
\widetilde{\G}^\ast = \DcNeg\gQ
\quad\mathrm{with\;its\;inertia\;}\quad
{\widetilde{\G}}^{\ast}{}\transpose\Dc\widetilde{\G}^{\ast} = \eye.
\end{align}
These asymmetric projections of CA/MCA/DiSCA/DiMCA 
are preserved in the gsGSVD with sparsification.

\begin{property}
Distributional equivalence: Rows (or columns) 
proportional to each other can be replaced
by their sum without affecting the results of the analysis.
\end{property}

Because these methods analyze the frequencies of the occurrences, 
two rows (or two columns) \ha{proportional 
to each other have identical profiles 
\citep[i.e.,
see][for details.]{escofier1965}}
Therefore
these rows (or columns) are represented by 
\textit{two} points having the same coordinates
in the component space, 
and so  these two points can be merged into \textit{one} 
whose mass is the sum of their original masses
\citep{escofier1965,fichet2009metrics,benzecri1973analyse,greenacre1984theory}.
In addition, merging these two points 
does not change the geometry of the component space.
This property is kept in gsGSVD with sparsification.

\begin{property}
Barycentric projection: 
the barycenters of the row factor scores, 
and the column factor scores are equal to the null vector 
\citep{escofier1965}.
\end{property}
The row and the column factor scores of
CA/MCA/DiSCA/DiMCA share a common barycenter 
(i.e., weighted mean) of 0; Formally
\begin{equation}\label{barycent_sCA}
\frac{1}{I} \rr\transpose \getc{\f}{\ell} 
= \frac{1}{J} \cc\transpose \getc{\g}{\ell} = 0.
\end{equation}
Specifically in MCA and DiMCA, 
the barycenter of each variable (i.e., set of columns) is null:
\jcy{\begin{equation}\label{barycent1_sMCA}
\sum_{j = 1}^{J_k} \getsc{c}{j}{k} \g_{j, k, \ell} = 0,
\end{equation}
where $\getsc{c}{j}{k}$ is the column weight for the $j$th level of the $k$th variable, 
and $\g_{j, k, \ell}$ is the factor scores of the $j$th level of the $k$th variable 
on the $\ell$th component.}
In addition, all observations that belong to a given variable level,
have their mean row factor scores equal to the column 
factor score of this variable level. 
\jcy{This property is lost in sCA and sDiSCA 
because the $\mathrm{proj}_{\mathcal{L}_1}$ operator 
sparsifies the variable levels individually on rows and/or columns.}
However, 
this property is kept in sMCA and sDiMCA, 
\jcy{because these methods use internally
the $\mathrm{proj}_{\mathcal{L_G}}$ 
operator to keep or sparsify 
the levels of each variable in groups, 
which conserves their barycenters. 
Simulation results are illustrated in \ha{*** Figure XX.***}}


\begin{property}
Embedded solution: 
The GSVD of the non-centered matrix (i.e., $ \Z $) 
has a first singular value of 1, 
a first left generalized singular vector equal to $ \rr $, 
and a first right generalized singular vector equal to $ \cc $. 
In addition, the subsequent 
components are the same
as those of the GSVD of the double centered matrix $ \X $.
\end{property}

The embedded solution holds because the GSVD in CA 
(Equation \ref{CA_GSVD}) can be rewritten as:
\begin{equation}\label{rewriteGSVD_sCA}
\begin{gathered}
\X = \Z-\rr\ct = \gP \BDelta \gQt = 
\sum_{\ell = 1}^L \getc{\delta}{\ell} 
\getc{\gp}{\ell} \getct{\gqt}{\ell} \\
\mathrm{under\;the\;constraints}\quad
\gPt\DrNeg\gP = \gQt\DcNeg\gQ = \eye
\end{gathered}
\end{equation}
which gives

\begin{equation}\label{noncentGSVD_sCA}
\begin{gathered}
\Z = \rr\ct + \gP\BDelta \gQt =  1 \times \rr\ct +  
\sum_{\ell = 2}^L \getc{\delta}{\ell} 
\getc{\gp}{\ell} \getct{\gqt}{\ell} \\
\mathrm{under\;the\;constraints}\quad
\gPt\DrNeg\gP = \gQt\DcNeg\gQ = \eye.
\end{gathered}
\end{equation}
Therefore, when the non-centered data $ \Z $ is analyzed, 
the first generalized singular value $\getc{\delta}{1}$ equals 1, 
the first left generalized singular vector $\getc{\gp}{1}$ equals $ \rr $,
and the first right generalized singular vector $\getc{\gq}{1}$ 
equals $ \cc $. With $ \rr\ct $ 
computing the \textit{expected} frequencies of $ \Z $ 
under independence, the CA of $ \X $, where $ \X = \Z- \rr\ct$, 
analyzes the deviation of the \textit{observed} data (i.e., $ \Z $) 
from independence (i.e., as given by matrix $ \rr\ct $).

When the singular vectors are sparsified, 
the embedded solution of the CA/MCA/DiSCA/DiMCA (Property
\ref{prop:embeddedsolution})
framework does not hold. 
Therefore, when the non-centered data $ \Z $ 
is analyzed by sCA/sMCA/sDiSCA/\-sDiMCA, 
the first pseudo-singular-value $\getc{\delta}{1}$
will be close to, but smaller than, 1, 
the first left (respectively right) 
generalized pseudo-singular-vector $\getc{\gp}{1}$ 
(respectively $\getc{\gq}{1}$)
will be close to but not equal to  $ \rr $ (respectively $\cc$). 

\subsection{Evaluating the sparsification}\label{section:SI}

The most straightforward way to evaluate the sparsity 
of left and right pseudo-generalized singular vectors 
is to express their number of zeros
as a function of $s$ (i.t., the sparsity parameter). 
\cite{Liu2023}
noted that such a crude index of sparsity 
would not be useful for choosing an ``optimal'' 
value for the sparsity parameters.
\jcy{Recently,} \citeauthor{Trendafilov2017}'s 
(\citeyear{Trendafilov2017})
sparsity index combines two measures: 
a measure of sparsity 
and a measure of how close the reduced 
rank sparse matrix is to the original data matrix. 

We measure sparsity with three different 
indices: 
1)
$\vartheta_{\P}$, 
the ratio of the number of zeros to the
total number of coefficients 
in the loading matrix $\P$, 
and 
2)
its counterpart $\vartheta_{\Q}$, 
the ratio of the number of zeros to the total number of coefficients
in the loading matrix $\Q$, and 
3)
$\vartheta$, 
the ratio of the number of zeros in both $\P$ 
and  $\Q$ to the total number of coefficients in $\P$ and  $\Q$.
\begin{align}
\nonumber
  \vartheta_{\P} & = \frac{\#_{0}(\P)}{I \times \ell} \\[2ex]
  \nonumber
  \vartheta_{\Q} & = \frac{\#_{0}(\Q)}{J \times \ell} \\[2ex]
  \vartheta & = \frac{\#_{0}(\P) + \#_{0}(\Q)}{(I+J) \times \ell}
\end{align}

\jcy{Given $L$ being the specified number of dimensions to estimate, }we measure the fit, noted $\hat\tau$, as 
the ratio of the sum of the $L$ squared pseudo-singular
values to the sum of the first $L$ squared singular values: 
\begin{equation} 
  \hat\tau = \frac{\displaystyle \sum_{\ell=1}^L \hat 
  \delta_\ell^2}{ \displaystyle \sum_{\ell=1}^L \delta_\ell^2}
\end{equation}

Combining both two types of ratios (sparsity and fit), 
gives three different sparsity indices:
\begin{align}
\nonumber
  \varsigma_{\P} & = \vartheta_{\P} \times \hat\tau \\\nonumber
  \varsigma_{\Q} & = \vartheta_{\Q} \times \hat\tau \\
  \varsigma & = \vartheta \times \hat \tau,
\end{align}
where $\varsigma_{\P}$ is the compromise between fit 
and sparsity for the left generalized singular vectors, 
$\varsigma_{\Q}$ is the compromise between fit and sparsity
for the right generalized singular vectors, 
and
$\varsigma$ is the compromise between fit 
and sparsity for both the left and 
the right generalized singular vectors.

Depending on the application, the analyst
can use one or more 
of
 these three sparsity indices to
 select the appropriate value 
 of the sparsity parameters $s_{\p}$ and  $s_{\q}$. 
 We illustrate these indices in the result section.

\begin{figure}[!h]
    \centering
    \includegraphics[width=.5\textwidth]{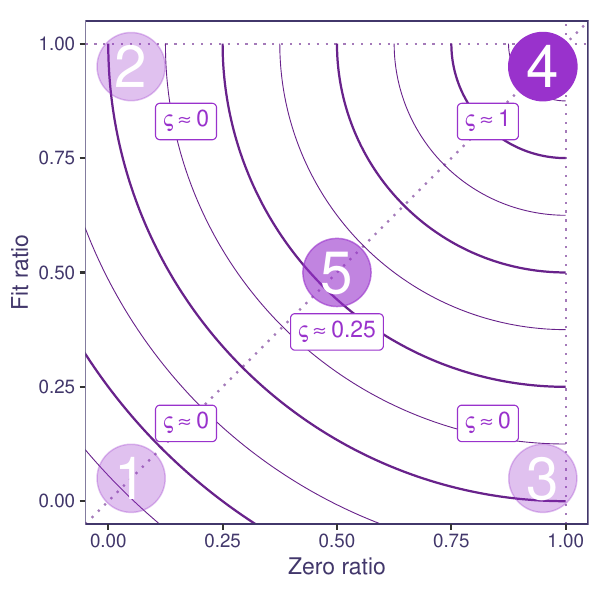}
    \caption{Graph representing different possible values 
    for the sparsity index on a map of the 
    ``fit'' as a function of the ``zero ratio.'' 
    The five zones represent five possible combinations 
    of the two ratios, 
    along with the corresponding values of the sparsity index $\varsigma$.}
    \label{fig:sparseIndex}
\end{figure}

Figure \ref{fig:sparseIndex} 
depicts the range 
of possible values for the sparsity indices 
on a graph representing the ratio of zeros
on the $x$-axis and the ``fit'' on the $y$-axis.
On this figure, the result of a sparse analysis 
would be represented as a dot, 
according to how sparse the loadings are (``zero ratio'') 
and how close the lower rank sparse decomposition
of the data is to the original data (``fit''). 
We split the graph into 5 zones: Zones 1 to 3 correspond 
to a low sparsity index because either or both ``fit'' 
and ``zero ratio'' are close to zero, 
Zone 4 corresponds to a sparsity index close to its 
maximum value of 1, 
when very few variables are selected and they represent 
most of the information in the data, and, finally,
Zone 5 corresponds to a middle ground,
where a compromise is reached between sparsity and fit.
In the following \jcy{Result section}, we will provide such map 
for each example. 
The optimal value for the sparsity parameters 
and the number of dimensions 
will be chosen by maximizing the sparsity index.


\section{Results}

\subsection{sCA}

\subsubsection{Data}

We applied sparse Correspondence Analysis 
to a contingency table
storing the number of deaths in the USA in 2018
as a function of age and causes of death. 
The causes of death were categorized 
according to the International Statistical Classification of Diseases 
and Related Health Problems (ICD) \citep{icd10}. 
We grouped the ages of death by every 5 years 
starting from the age of 1 and summing all deaths 
above 100 years old into a single group called ``100+.'' 
In addition, we removed the deaths 
before the age of 1 because \jcy{they mostly belong to} a single category of the perinatal causes. 
\jcy{Including this age group would lead to this specific cause-age association dominating the result and shadowing other effects as} deaths resulting from perinatal 
causes were rare in other age groups 
(\jcy{i.e., 57 out of 118 cases were before 1 year old with the other 61 cases spread across the other 21 age groups.}). 
The resulting data analyzed by CA and sCA was a 21 
(age groups) by 19 (causes of death) 
contingency table with counts indicating the number 
of deaths of each cause at each age range.

\subsubsection{Results}

\begin{figure}[!!h]
    \centering
    \includegraphics[width=\textwidth]{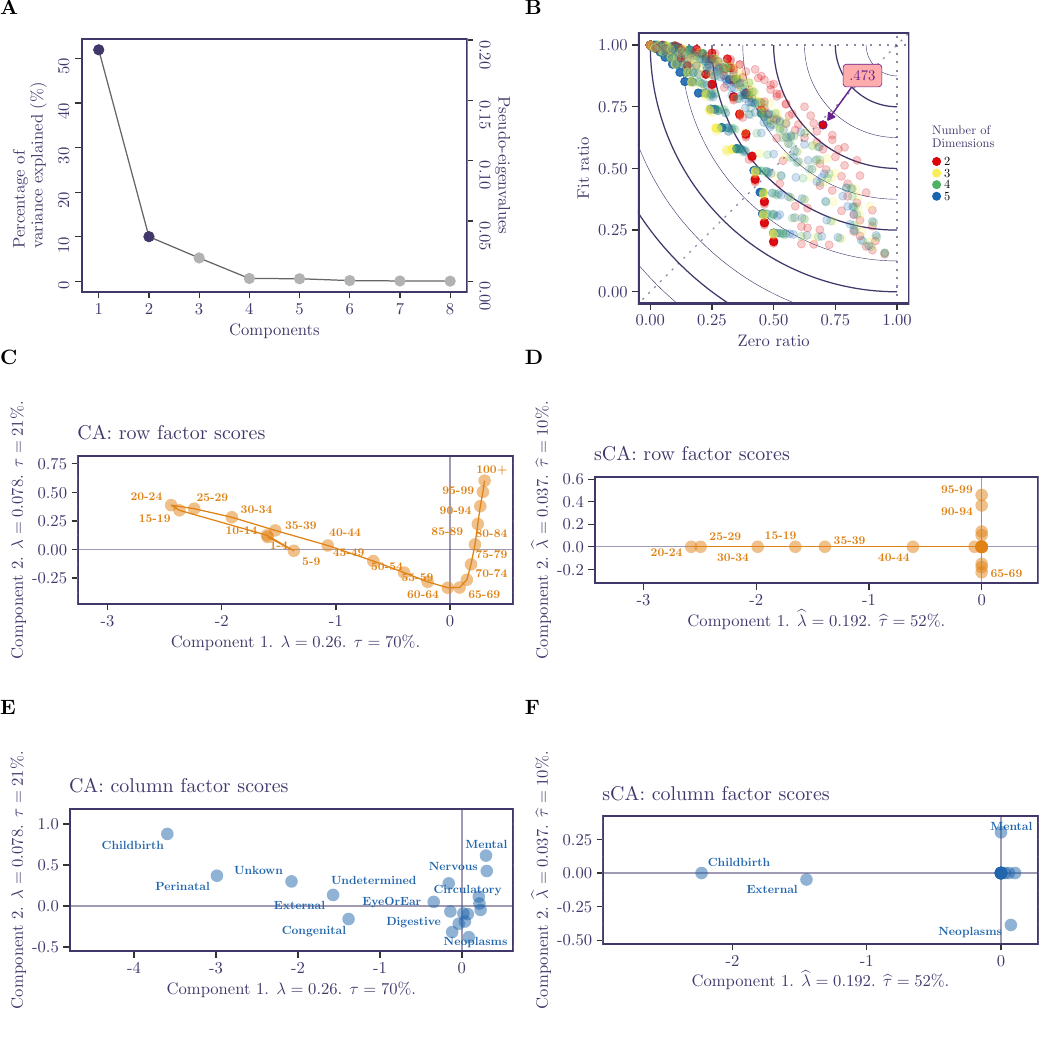}
    \caption{Results of CA and sCA. 
    (A) shows the scree plot of sCA with the optimal 
    number of components colored in purple. 
    (B) shows the fit-to-zero-ratio plot 
    and highlights the optimal solution 
    that has the maximum sparsity index. 
    (C) and (D) show the row factor scores 
    (which represent different age ranges) 
    from CA and sCA. 
    (E) and (F) show the column factor scores 
    (which represent different causes of death) from CA and sCA.}
    \label{sCAres}
\end{figure}

Figure \ref{sCAres}A shows the scree plot of the sCA result.
Here, we only consider solutions 
\ha{with}
\jcy{more than} 2 components 
\jcy{and found that sCA with 2 components gave the optimal result}, 
\jcy{which has}
the sparsest solution with the largest fit.
\jcy{These characteristics are indicated in 
Figure \ref{sCAres}B, which illustrates 
the sparsity index of possible solutions 
with the chosen one being }
closest to the upper right corner 
where both fit and zero ratios are equal to 1. 
The optimal sparsity parameter for the age groups is 
$ .51 \times \sqrt{21}$ (21 age groups), 
and the optimal sparsity parameter for the causes of death is 
$ .31 \times \sqrt{19}$ (19 causes of death). 
The sparsity index from this analysis was equal to $.473$.

\clearpage

In CA, the first component has an eigenvalue 
of .26 which explains 70\% of the inertia, 
and the second component has an eigenvalue of .08 
which explains 21\% of the variance. 
The row and the column factor scores show 
that the first component is characterized 
by deaths between the late teens and twenties 
that are related to childbirth and 
deaths at a young age (< 35 years old) that are related 
to congenital and perinatal conditions 
(see the horizontal axis in Figures \ref{sCAres}C and E). 
The second component is characterized by deaths 
at older ages (> 65 years old)
that are mostly driven by malfunction 
in the nervous system 
such as Alzheimer's and Parkinson's diseases 
(see the vertical axis in Figures \ref{sCAres}C and E). 
With sCA, the components give a clearer pattern 
that facilitates interpretation. 
The sparsified row and the column factor scores show 
that the first component is characterized by deaths
in early adulthood (between 20--35) that are related
to childbirth and external causes (e.g., suicide). 
The second component \jcy{differentiates} deaths 
between 60--80 years old, 
which relate to the causes of neoplasms 
(e.g., cancer), 
from deaths at 90 years old, 
which relate to mental conditions. 

\subsection{sDiSCA}
\subsubsection{Data}
To demonstrate sparse DiSCA, 
we used a data set that counts the numbers of punctuation marks 
in works of authors from different times and origins. 
The texts were extracted from the Gutenberg Project 
using the R-\texttt{gutenbergr}
package \citep{robinson2021gutenberg}. 
In this data set, we included authors from France 
($N = 63$), the United Kingdom ($N = 61$), 
and the United States ($N = 36$) 
and from different periods (see Table 1).
The data counted the number of occurrences 
of the following nine
punctuation marks: comma (,), period (.), 
question mark (?), exclamation mark (!), 
colon (:), semicolon (;) apostrophe ($^\prime$), 
quotation marks for both English (`` ") and French (« »),
dashes (-), and M-dashes (---).
We did not include the translated works and 
only considered the books written in the authors' original languages. 
We used DiSCA and sDiSCA to perform discriminant 
analysis on the eight author groups listed in Table \ref{tab:author}.

\begin{table}[!h]
\centering
\caption{Numbers of authors in each group.}
\begin{tabularx}{0.5\textwidth}{c *{3}{Y}}
\toprule
Origin & \multicolumn{3}{c}{Time Period (centry)} \\
& 18th and before & 19th & 20th and after \\
\midrule
France & 16 & 28 & 19 \\
UK & 11 & 28 & 22 \\
\cmidrule{2-3}
US & \multicolumn{2}{c}{14} & 22 \\
\bottomrule
\end{tabularx}

\bigskip
\small\textit{Note}. UK: United Kingdom; US: United States.
\label{tab:author}
\end{table}

\subsubsection{Results}
\begin{figure}[!!h]
    \centering
    \includegraphics[width=.94\textwidth]{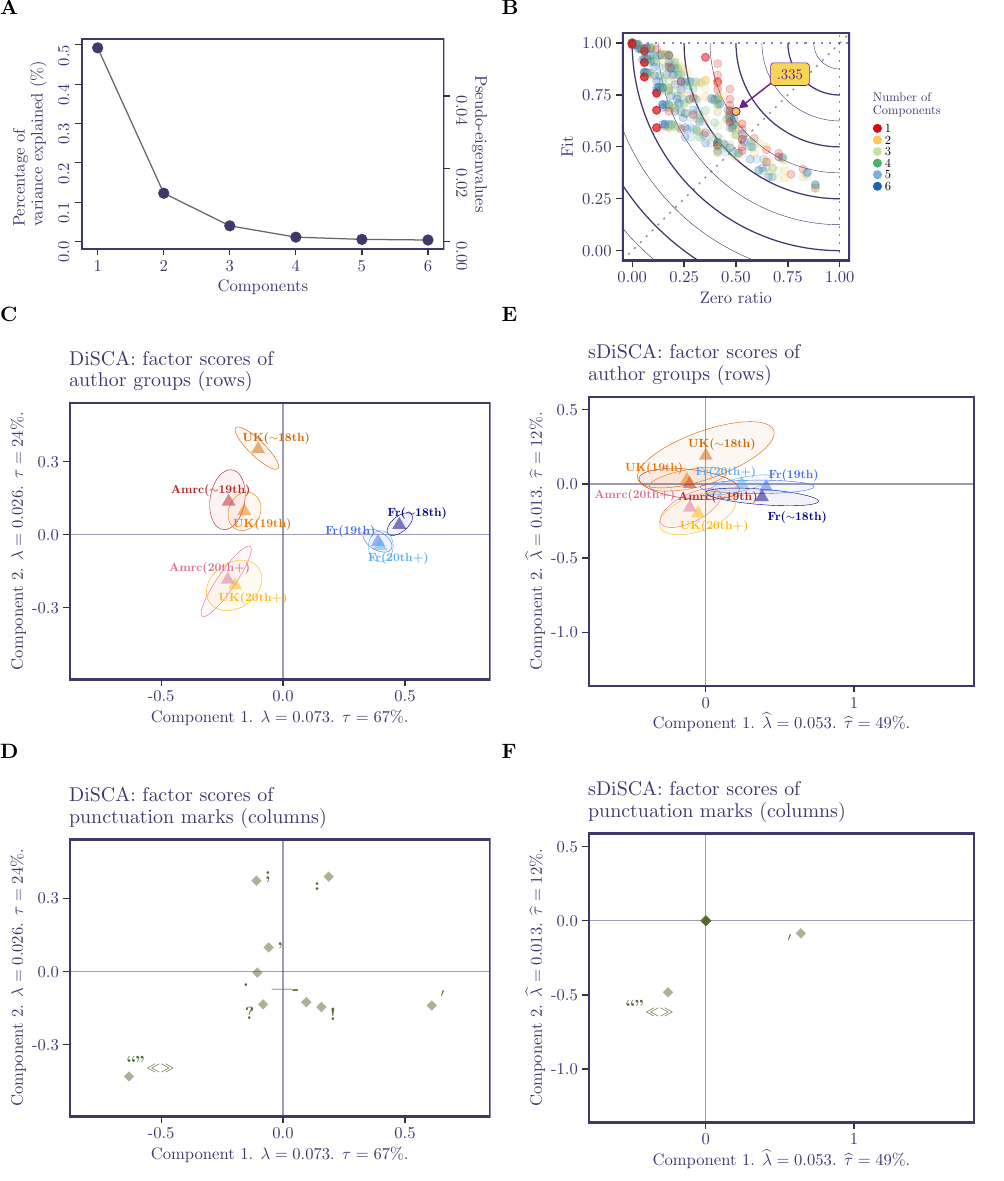}
    \caption{Results of DiSCA and sDiSCA. 
    (A) shows the scree plot of sDiSCA with the optimal number of components. 
    (B) shows the fit-to-zero-ratio plot 
    and highlighted the optimal solution with the maximum sparsity index. 
    (C) and 
    (D) show the row factor scores 
    (which represent the groups of authors) from DiSCA and sDiSCA. 
    (E) and (F) show the column factor scores
    (which represent the punctuation marks) from DiSCA and sDiSCA.}
    \label{sDisCAres}
\end{figure}

The sDiSCA results are shown in Figure \ref{sDisCAres}. 
Here, we consider a solution with 2 components because
\jcy{it} gives \jcy{the}
optimal results according to the sparsity index 
(Figure \ref{sDisCAres}B). 
As shown in Figure \ref{sDisCAres}B, 
this result is also the sparsest solution with the largest fit 
(closest to the upper right corner where both fit and zero ratios equal 1).
The sparsity parameter for the individuals is set to equal $.71 \times \sqrt{8}$ 
(8 author groups), 
and the optimal sparsity parameter 
for the punctuation marks is $ .41 \times \sqrt{9}$ (9 punctuation marks). 
The sparsity index from this analysis is equal to $.335$.

The original and the sparsified component spaces are shown 
in Figure \ref{sDisCAres}C--\ref{sDisCAres}F. 
Figure \ref{sDisCAres}C shows the factor scores of the author groups 
from DiSCA, and Figure \ref{sDisCAres}D shows 
the same results from sDiSCA. 
The individual authors are projected 
as supplementary observations onto the same components 
and are used to illustrate the stability of these \jcy{group means}
by estimating their 95\% bootstrap confidence intervals. 
When two confidence intervals overlap,
the difference between the two groups is considered non-significant
(with $ \alpha = .05$). 
In DiSCA, the first component distinguishes the French-speaking 
authors from the English-speaking authors,
and the second component distinguishes the 
English-speaking authors from different periods of time. 
By contrast, 
the French authors are more consistent across time. 

The results from sDiSCA show a similar pattern 
in Figure \ref{sDisCAres}D; 
the first component identifies 
only the three French groups of authors, 
and the second component identifies the oldest 
and the most recent groups of English authors. 
However, the bootstrap confidence 
intervals have quite a large variance---%
mostly because these intervals are derived 
from the supplementary projections of the authors
with sparsification involved. 
While the sparsification algorithm sparsifies the number 
of contributing groups of each component\jcy{,} 
the group separations are, as a trade-off, no longer optimized, 
and \jcy{thus results in supplementary projections
are less segregated between groups.} 

The same trade-off is also reflected in the classification accuracy. 
In DiSCA, the discrimination between groups 
is measured by the accuracy of correct classifications. 
To classify observations, their supplementary factor 
scores are compared to the group factor scores 
and classified as belonging to the closest group. 
The accuracy rate of DiSCA is .44 
and the accuracy rate of sDiSCA is .16. 
Although  lower than for  DiSCA, 
the accuracy is still above chance level (.125). 
The accuracy rate is higher with respect to the origins 
(DiSCA: .74; sDiSCA: .47; chance level: .33) 
with the French authors being classified most accurately 
(DiSCA: 1; sDiSCA: .68) as compared to the UK (DiSCA: .59; sDiSCA: .52) 
and the US (DiSCA: .56; sDiSCA: .27) authors.

Figure \ref{sDisCAres}E shows the factor scores of the punctuation marks 
from plain DiSCA, and Figure \ref{sDisCAres}F shows the same results from sDiSCA. 
The language effect on the first component is associated 
with the difference between how these authors punctuated 
quotation marks and apostrophes. 
This language effect is therefore consistent 
with the difference in how apostrophes are used differently 
in English (to represent possession) and French (for grammatical purposes). 
The time effect on the second component is associated with 
how authors used colons, semicolons, and commas 
(which connect sentences) versus
how they used quotation marks, apostrophes,
question, and exclamation marks (which end sentences). 
Therefore, the time effect could be related to how the style changes 
in English writing across time. 
Similar punctuation patterns are also shown 
in the results of sDiSCA which identified  
the quotation marks and the apostrophes as
 contributing most to the group differences.

\clearpage

\subsection{sMCA}
\subsubsection{Data}

To illustrate sMCA, we used a data set on the Chinese 
version of the 
Independent and Interdependent Self Scale (C-IISS) 
developed by \cite{lu2007developing}. 
This data set describes 130
undergraduate students 
(77 females and 53 Males;
$ M\mathrm{_{age}}$ = 19.49, and $sd\mathrm{_{age}}$ = 1.52) 
from National Cheng Kung University.
These participants signed written informed consents 
and received NTD 120 at the end of the experiment. 
The C-IISS comprised forty-two 7-point Likert scale questions 
(1 = strongly disagree; 7 = strongly agree). 
Among the 42 items, 21 of them measured independence 
(i.e., if one is aware of and values oneself as an individual) 
and the other 21 measured interdependence 
(i.e., how one values oneself and acts based on 
\jcy{one's} cohort). 
Before we analyzed the data with sMCA, because the results of MCA
(and sMCA) will be driven by the rarity of an event,
we binned the responses of each item into categories
of comparable sizes (see Supplementary Figure \ref{sMCABin}). 
The association patterns between these items are given in 
Supplementary Figure \ref{sMCAheat}. 
We analyzed the binned data with 
\jcy{both MCA and} sMCA\jcy{, where we only sparsified} the items. 
Although sparsifying the individuals 
(or both) is also possible, we did not sparsify 
them because the separation of individuals was 
not of interest in this study.

\subsubsection{Results}
\begin{figure}[!!h]
    \centering
    \includegraphics[width=0.93\textwidth]{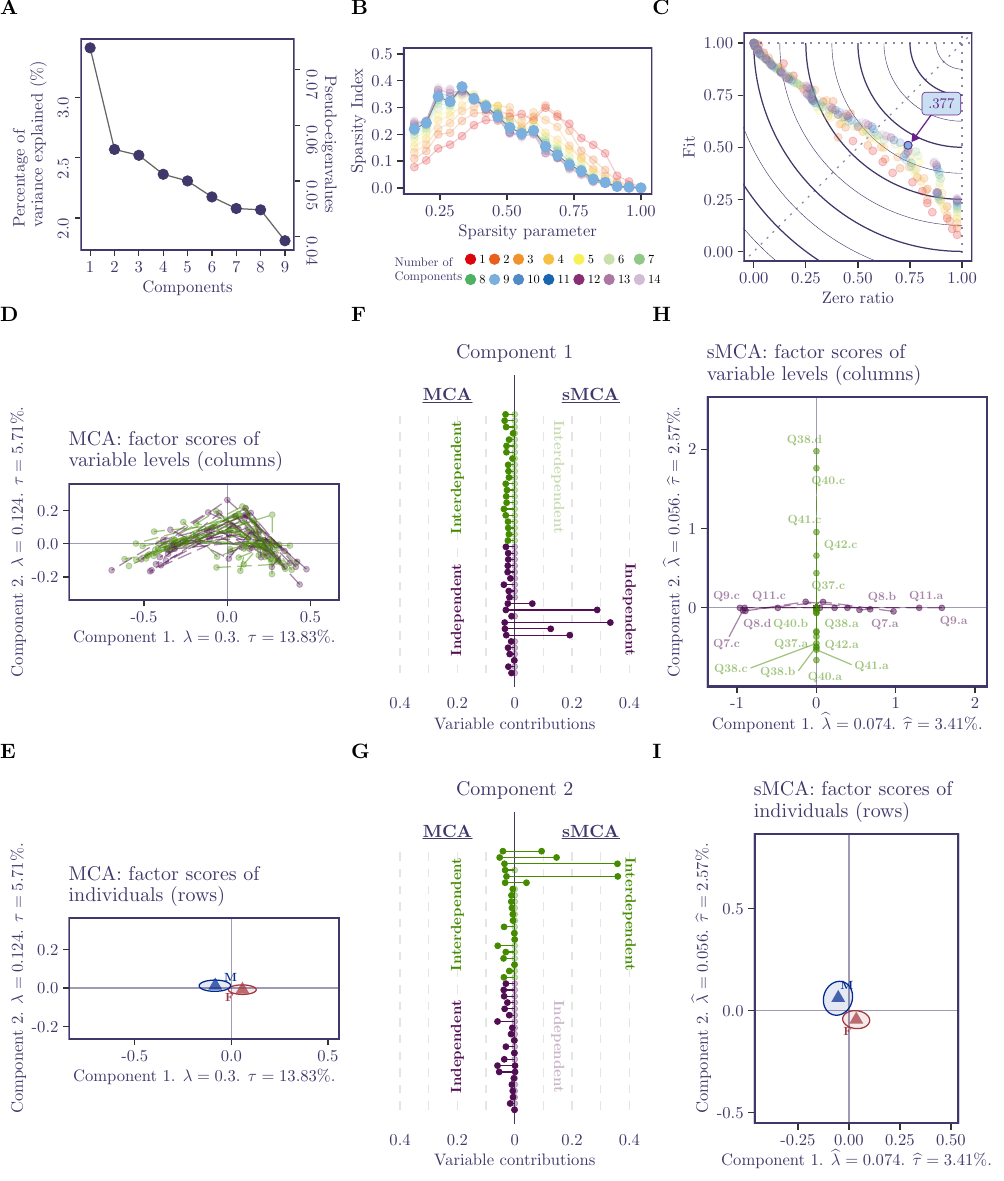}
    \caption{Results from MCA and sMCA. 
    (A) shows the scree plot of sMCA with the optimal number of components. 
    (B) and (C) show the fit-to-zero-ratio plot 
    and highlights the optimal solution with the maximum sparsity index. 
    (D) and 
    (H) show the column factor scores 
    (which represent variable levels) from MCA and sMCA 
    with the levels of the same variables connected 
    from low to high to illustrate the trend. 
    Items measuring interdependent self 
    are colored in green, and items 
    measurng independent 
    self are colored in purple. 
    (F) and 
    (G) show the variables that contribute 
    to the two components of MCA and sMCA side by side. 
    (E) and (I) 
    show the mean row factor scores (which represent individuals) 
    from MCA and sMCA with the blue triangle denoting 
    the male group and the red triangle denoting the female group. 
    The ellipses around the means illustrate the 
    95\% bootstrap confidence intervals of the group means.}
    \label{sMCAres}
\end{figure}

The sMCA results are shown in Figure \ref{sMCAres}.
For this data set, sMCA with 9 components 
gives the optimal results according to the sparsity index 
(Figure \ref{sMCAres}B). 
As shown in Figure \ref{sMCAres}C, 
this result is also the sparsest solution
with the largest fit (closest to the upper 
right corner where both fit and zero ratios equal 1). 
The sparsity parameter for the individuals was set to
$\sqrt{130}$ (i.e., no sparsity; 130 individuals), 
and the optimal sparsity parameter for the items is $ .33 \times \sqrt{42}$
(42 item levels). The sparsity index from this analysis equals $.377$.

The plain and sparsified component spaces are shown in 
Figure \ref{sMCAres}D-\ref{sMCAres}I.
In plain MCA---with 
Benzécri's
\citeyearpar{benzecri1979calcul}
eigenvalue correction---%
the first component eigenvalue of .08 
explains 72\% of the inertia, and 
the second component with an eigenvalue of .01 
explains 9\% of the inertia. 
Because it is still unclear how such a correction 
should be applied to sMCA,
we report here the uncorrected values directly from the 
GSVD for easier comparison between the two methods. 
Without Benz{\'e}cri's correction, 
the first component from plain MCA with 
an eigenvalue of .30 explains 13.83\% of 
the variance, 
and the second component with an eigenvalue of .12 explains 5.71\% of 
the variance. 
The factor scores of the individuals 
are grouped according to their sex at 
birth, and the factor scores of the items 
are colored according to their corresponding 
category. 
The results showed that MCA generated components 
that could be difficult to 
explain given the complex pattern of loadings 
(see Figure \ref{sMCAres}D), 
whereas the 
first two components from sMCA 
distinguished items from different 
categories 
(see Figure \ref{sMCAres}F-\ref{sMCAres}H). 
Specifically, the 
first component of sMCA identifies 
linear level effects of 5 questions that 
measure the level of independence; 
the second component of sMCA identifies linear 
level effects of 5 questions 
that measure the level of interdependence. 
The first 
component from sMCA with an eigenvalue of .07 
explains 3.41\% of the variance, 
and the second component with an eigenvalue 
of .06 explains 2.57\% of the 
variance. 
Figures \ref{sMCAres}E and  \ref{sMCAres}I 
show the 
factor scores of the individuals 
with the mean factor score of each sex group. 
To examine group effect post-hoc, 
the stability of these group means is illustrated by their 95\% bootstrap 
confidence intervals; when two confidence intervals overlap, the difference 
between the two groups is considered non-significant (with $ \alpha = .05$) \jcy{in MCA but significant in sMCA}. 
\jcy{In general,} the individual factor scores from both analyses showed similar group separations \jcy{with a marginal difference}. 

\clearpage

\subsection{sDiMCA}

\subsubsection{Data}

\jcy{The data we used to demonstrate DiMCA and sDiMCA} 
were collected from undergraduate students enrolled 
in an introductory psychology statistics course 
taught by various instructors, 
at an urban public college in the northeast of the United States. 
The course covered descriptive statistics, 
hypothesis testing, 
and an overview of advanced statistical procedures. 
The Math Assessment for College Students 
\citep[MACS;][]{Rabin2018}
was administered to 460 participants across five semesters, 
measuring basic mathematics skills through 
a 30-item paper-and-pencil test, 
which covered five general content domains 
(Table~\ref{tab:math}). 
Informed consent was obtained under 
an IRB-approved protocol and participants were not compensated. 
Demographic and academic performance data were also collected. 
The goal of the study was to examine 
the relationship between basic mathematics skills, 
demographic data, and academic performance.

The MACS questions took approximately 
40 minutes and students completed 
the MACS during the first week of the semester. 
All 30 MACS items were graded
with no partial credit, with each response 
recorded as 0 = incorrect and 1 = correct 
\citep[scored by a single rater and re-scored 
by a second independent rater;][]{Rabin2018}. 
The course was computationally based,
with students learning to perform 
statistical tests manually and use 
statistical software programs. 
Academic performance was evaluated 
based on the average score computed 
from three exams conducted in each semester. 
All exams included multiple-choice 
and problem-solving questions 
that covered basic statistics theory and applications.
The average score was then categorized 
into a letter grade common to most
undergraduate statistics courses, 
where: A = 90-–100\%; B = 80--89\%; C = 70--79\%; 
D = 60–-69\%; and F = below 60\%.
For this illustration, 
A and B grades were grouped together 
as were D and F grades, resulting 
in a total of 3 groups (i.e., AB, C, DF). 

\begin{table}
\centering
  \caption{Description of the groups for the Math questionnaire data}
  \label{tab:math}
  \smallskip
  \begin{tabular}{llr}
    \toprule
    Domain & Short Name & Number of Items \\
    \midrule
    basic algebraic skills & Algebra & 8 \\
    basic arithmetic skills & Arith & 5 \\
    categorization and ranges & CatRang & 3 \\
    decimals, fractions, and percentages & DecFracPerc & 12 \\
    visual understanding & Visual & 2 \\
    \bottomrule
  \end{tabular}
  \smallskip
\end{table}

\subsubsection{Results}

\begin{figure}[!!h]
    \centering
    \includegraphics[width=\textwidth]{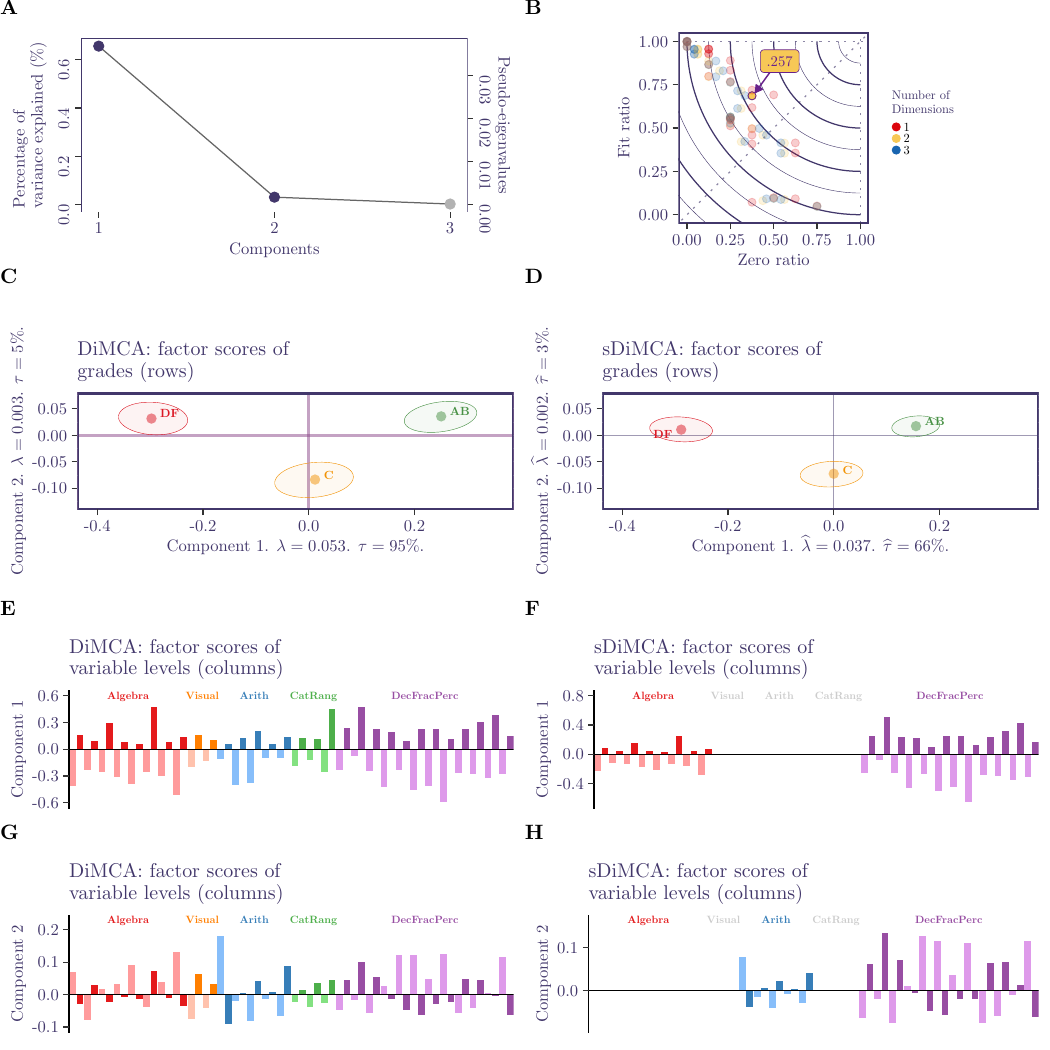}
    \caption{Results from DiMCA and sDiMCA. 
    (A) scree plot of DiMCA with components to keep (dark purple). 
    (B) fit-to-zero-ratio plot with optimal solution highlighted. 
    (C) and (D) row factor scores from (respectively) DiMCA and sDiMCA. 
    Ellipsoids show bootstrap derived confidence intervals. 
    (E) to (H) factor scores for the variables 
    colored by category for DiMCA (on the left) 
    and sDiMCA (on the right)---only non-null values are shown.}
    \label{sDiMCAres}
\end{figure}

DiMCA and sDiMCA results
are shown in Figure~\ref{sDiMCAres}. 
For this data set, sMCA with 2 components 
gives the optimal non-unidimensional 
results according to the sparsity 
index (Figure \ref{sDiMCAres}B). 
As shown in Figure \ref{sDiMCAres}B, 
this result is also the sparsest 
solution with the largest fit 
(closest to the upper right corner where both fit and zero ratios equal 1). 
The sparsity parameter for the groups 
is set to equal $.79 \times \sqrt{5}$, 
and the optimal sparsity parameter 
for the items is $ .57 \times \sqrt{60}$ (60 item levels). 
The sparsity index from this analysis equals .257.

\clearpage

As shown in Figure~\ref{sDiMCAres}, 
both DiMCA and sDiMCA differentiated between 
AB-level and DF-level students along 
the first dimension, and this difference 
was driven by the overall number 
of correct and incorrect responses 
on the MACS, where AB-level students 
had more correct responses and DF-level students 
had more incorrect responses. 
\jcy{Specifically}, 
DiMCA indicates that responses 
for all 30 MACS items reliably 
contribute to differences in performance, a
pattern
implying that DF-level students---%
compared to AB-level students---%
had an overall low performance across 
all basic domains of mathematics 
In contrast, sDiMCA reveals that differences 
in performance on particular items related 
to basic algebraic skills and decimals, 
fractions, and percentages most discriminated
between AB-level students and DF-level students. 
These topics are usually covered 
in elementary and middle school and may not have been reviewed 
prior to college-level statistics.
Therefore, students could be encouraged 
to take statistics earlier than later 
in college when basic mathematics skills 
are more likely to be remembered. 

Based on MACS performance, 
certain students could be identified as ``under-prepared'' 
or ``at risk'' for poor outcomes 
in an introductory statistics course 
and might benefit from remediation 
before attempting courses that rely heavily on mathematics. 
Such remediation could be developed 
either as in-person or online formats 
with assignments completed independently 
or with a small learning group. 
Ideally, such efforts would provide 
opportunities for early success 
in introductory statistics courses through mastery 
of relevant pre-requisite mathematics skills. 

In addition, both DiMCA and sDiMCA differentiated 
the C-level students from other grade levels 
along the second dimension, but DiMCA, 
once again, showed that responses for 
all 30 MACS items reliably contributed to this difference. 
However, sDiMCA revealed that C-level students 
were able to correctly answer specific 
subsets of items in basic arithmetic 
(e.g., order of operations) and decimals, 
fractions, and percentages (e.g., adding fractions) 
that differentiated them from DF-level 
students, but incorrectly answered enough items 
(e.g., rounding to the nearest decimal, decimal 
and percentage conversions) 
to differentiate them from the AB-level students. 
These results suggest that while C-level 
students are slightly better than 
DF-level students on certain concepts, 
these students might also benefit
from remediation so that all students
have the opportunity to begin the course 
with a similar level of mathematics proficiency.

\section{Conclusion and discussion}

In this paper,
we extended the sparse GSVD introduced in \citet{Yu2023}
to create the group-sparse GSVD (gsSVD) algorithm 
which keeps the orthogonality constraints 
of the GSVD \jcy{while integrating}
sparsification \jcy{with} metric and group constraints.
We applied the gsGSVD to sparsify CA, MCA, 
and their respective discriminant analysis versions: DiSCA and DiMCA. 
We illustrated
these sparsification methods
and compared their analytical merits on four real data sets.
We also integrated group constraints into the 
sparsification problem---%
an essential property for categorical data analysis where
variables are represented by blocks of columns.
These group constraints could also 
be applied to different situations
such as 
for data with
\emph{a priori} groups of variables.

This new gsGSVD algorithm seeks optimum solutions that satisfy strong 
constraints such as orthogonality and sparsity.
\jcy{Interestingly, our results also
showed that sparsification can be interpreted 
as a form of cluster analysis
performed on the rows or the columns 
of the data matrix---an approach akin
to, for example,  
spectral clustering which partitions the variables 
based on their singular vectors \citep{Kannan2000OnCB}
or graph analysis \citep{spielman2011spectral}.}
\jcy{Overall, we demonstrated how 
this new algorithm provides new multivariate tools 
to explore the complex structures of qualitative data.} 

\jcy{To sparsify solutions while maintaining orthogonality 
between components, 
we used POCS to integrate two sets of strong constraints: 
one for orthogonality and one for sparsification.
However, in some cases, these constraints can be so strong that 
they may not be all satisfied together.
}
For example, to sparsify the solution while keeping 
the components orthogonal,
we project the data onto the orthogonality constraints 
first and the sparsification constraints last. 
This approach guarantees that sparsity 
and orthogonality are satisfied to a high level of precision 
i.e., (below 1e-10). 
In practice, it is also possible 
to switch the order of projections 
and project onto the orthogonality constraints last. 
This order of projections prioritizes orthogonality 
over sparsity. It is worth noting that, sometimes, 
the solution of prioritizing different constraints 
could lead to different results. 
In certain circumstances, meeting both orthogonality 
and sparsity constraints is not feasible. 
In such cases, the prioritized constraint 
will be satisfied at the expense of the other. 

\jcy{It is worth noting} that the dimensions derived from the sparse SVD 
algorithm might not always be obtained in decreasing order
of explained variance---
a problem already present in \cite{witten2009penalized}. 
The original order should be kept 
to implement the orthogonality constraints and
for the transition formulas to work correctly. 
However, for convenience and to follow traditions,
the current output of the R-implementation
reorders the dimensions according to their 
variance. 

\vg{
\jcy{Additionally}, as is done by \jcy{\citep{witten2009penalized}} 
and \citep{guillemot2019constrained}, 
to make \jcy{Algorithm \ref{algo:gsGSVD}} 
a convex optimization problem with convex constraints, 
the constraints on the (group-)$\Lyi$- and $\Ler$-norms \jcy{used} 
are inequality constraints---%
an interpretation
that could hinder the interpretability of the results\jcy{.
For example, when} the $\Ler$-norms of some loadings 
\jcy{are less than 1}. 
In practice, however, these constraints are saturated, 
meaning that \jcy{it is rarely the case} 
that the $\Ler$-norms of the resulting 
pseudo generalized singular vectors 
\jcy{are not exactly 1}, 
or that their $\Lyi$-norms 
\jcy{are not exactly the} given sparsity parameters. 
\jcy{For cases when the constraints are unsaturated 
(i.e., they are not compatible), 
the solutions will not be unique and could be unstable. 
}}

The proposed sparse methods still require future developments 
on \jcy{solutions to preserve more CA-MCA-related 
properties and on} inference analysis 
to evaluate the stability and validity of the sparse solution. 
The current work also opens avenues to explore the properties
of the POCS projecting operators 
such as the projection of supplementary elements, which currently 
uses
approximations. 


Overall, this study presents 
an exciting foundation for refining and expanding 
the use of sparsification techniques. 
Future avenues could extend the current methods to
 two data table analyses such as
 Partial Least Squares Correlation (PLSC) and its CA-like extension
 PLS-CA and PLS-MCA. 
 Future directions could also incorporate the consideration 
 of hierarchical structures of variables or observations 
 such as overlapping groups of grouped variables, 
 exemplified in Single Nucleotide Polymorphism (SNP) data structured into pathways.



\section{Acknowledgements}
JCY receives funding from the 
Discovery Fund postdoctoral fellowship 
of the  Centre of Addiction and Mental Health

LR and AK received funding from the PSC-CUNY TRADA Award (62058-00 50)

LR received funding from NSF Research Experiences for Undergraduates (REU) Award (2050755)

\bibliography{sparse_mca}

\newpage

\appendix

\section{CA and inertia}\label{append:inertia}

Recall that, for an $ I \times J $ contingency table,
the independence $ \chi^2 $ statistic is computed as:
\begin{equation}\label{chi2_CA}
\chi^2 = \sum\limits_{i,j} 
\frac{(\mathrm{Observed}_{i,j} - 
\mathrm{Expected}_{i,j})^2}{\mathrm{Expected}_{i,j}}
\end{equation}
or, in  matrix notation, the $ \chi^2 $ 
is associated to the matrix $ \X $ 
(with the notations from Equations \ref{probmat_CA} 
and \ref{doubelcent_CA}) and computed as:
\begin{equation}\label{chi2mat_CA}
\frac{1}{N}\chi^2 = \trace{ \DcNeghalf\Xt\DrNeg\X\DcNeghalf } 
= \trace{ \DrNeghalf\X\DcNeg\Xt\DrNeghalf }.
\end{equation}
From the GSVD, the inertia computed 
from the row factors $ \F $ gives $ \frac{1}{N}\chi^2 $:
\begin{equation}\label{inertiaF_CA}
\begin{split}
&\trace{ \Ft\Dr\F } \\
&= \trace{ \gQt\DcNeg\Xt\DrNeg\Dr\DrNeg\gP\BDelta } \\ 
&= \trace{ \DcNeg\Xt\DrNeg\gP\BDelta\gQt } \\
&= \trace{ \DcNeg\Xt\DrNeg\X } \\
&= \trace{ \DcNeghalf\Xt\DrNeg\X\DcNeghalf } \\
&= \frac{1}{N}\chi^2,
\end{split}
\end{equation}
and so does the inertia computed from the column factors $ \G $:
\begin{equation}\label{inertiaG_CA}
\begin{split}
&\trace{ \Gt\Dc\G } \\
&= \trace{ \gPt\DrNeg\X\DcNeg\Dc\DcNeg\gQ\BDelta } \\ 
&= \trace{ \DrNeg\X\DcNeg\gQ\BDelta\gPt } \\
&= \trace{ \DrNeg\X\DcNeg\Xt } \\
&= \trace{ \DrNeghalf\X\DcNeg\Xt\DrNeghalf } \\
&= \frac{1}{N}\chi^2.
\end{split}
\end{equation}

\newpage

\section{The SVD and the generalized SVD (GSVD)}\label{append:gsvd}

The SVD decomposes the data matrix $\X$ into three matrices:
\begin{equation}
    \X = \P\BDelta\Qt
    \quad\suchthat\quad
    \Pt\P = \Qt\Q = \eye,
\end{equation}
where $\P$ (respectively $\Q$) is the 
$I \times L$ (respectively $ J \times L $) 
matrix of the left (respectively right) 
singular vectors, and 
\BDelta is a diagonal matrix with singular values 
$\delta$s stored on its diagonal.
The SVD solves the following maximization problem:
\begin{equation}\label{Maximization_SVD}
\begin{split}
&\argmax_{\getc{\p}{\ell},
\getc{\q}{\ell}} 
\left( \getc{\delta}{\ell} = 
\getct{\pt}{\ell} \X\getc{\q}{\ell} \right)
\quad\st\quad\\[1ex]
&\begin{cases}
\getct{\pt}{\ell} \getc{\p}{\ell} &= 1, \\
\getct{\qt}{\ell} \getc{\q}{\ell} &=1,
\end{cases}
\quad\mathrm{and,\;for\;any\;}\ell  \neq\ell^\prime,\quad
\begin{cases}
\getct{\pt}{\ell} \getc{\p}{\ell^\prime} &= 0. \\
\getct{\qt}{\ell} \getc{\q}{\ell^\prime} &=0.
\end{cases}
\end{split}
\end{equation}
Here, $\getc{\delta}{\ell}$ is the
$\ell$th singular value and is associated 
to the $\ell$th left (respectively right) 
singular vector $\getc{\p}{\ell}$ and $\getc{\q}{\ell}$.

Similar to the SVD, the GSVD also decomposes $\X$ 
into three matrices but with row and column 
metric matrices included in the constraints:
\begin{equation}
    \X = \gP\BDelta\gQt
    \quad\suchthat\quad
    \gPt\M\gP = \gQt\W\gQ = \eye,
\end{equation}
where $\gP$ (respectively $\gQ$) 
is the matrix of the left (respectively right) 
generalized singular vectors, $\M$ 
(respectively $\W$) is the row 
(respectively column) metric matrix 
represented by a positive definite matrix.
The GSVD solves the following maximization problem:
\begin{equation}\label{MaximizationEll_GSVD}
\begin{aligned}
&\argmax_{\getc{\p}{\ell},\getc{\q}{\ell}} 
\left( \getc{\delta}{\ell} = \getct{\gpt}{\ell} \X \getc{\gq}{\ell} \right)
\quad\st\quad\\[1ex]
\quad\st\quad
&\begin{cases}
\getct{\gpt}{\ell} \M\getc{\gp}{\ell} &= 1, \\
\getct{\gqt}{\ell} \W\getc{\gq}{\ell} &=1,
\end{cases}
\quad\mathrm{and,\;for\;any\;}\ell\neq\ell^\prime,\quad
\begin{cases}
\getct{\gpt}{\ell} \M\getc{\gp}{\ell^\prime} &= 0. \\
\getct{\gqt}{\ell} \W\getc{\gq}{\ell^\prime} &=0.
\end{cases}
\end{aligned}
\end{equation}
The maximization problem of the GSVD 
is equivalent to the following SVD of 
the weighted $ \X $ (denoted by $ \widetilde{\X} $), where 
\begin{equation}\label{weightX_gsvd}
\begin{split}
    \widetilde{\X} = \Mhalf\X\Whalf = \P\BDelta\Qt \\
    \suchthat\quad\Pt\P = \Qt\Q = \eye,
\end{split}
\end{equation}
with
\begin{equation}\label{UV_gsvd}
\begin{cases}
  \gP &= \MNeghalf\P \\
  \gQ &= \WNeghalf\Q
\end{cases}
\quad\suchthat\quad
  \gPt\M\gP = \gQt\W\gQ = \eye.
\end{equation}

\newpage

\section{Properties of CA, MCA, DiSCA, and DiMCA}\label{append:properties}

Because of the specific preprocessing 
steps and the metric constraints, CA 
(and therefore MCA, DiSCA, and DiMCA) 
has several specific properties \citep{escofier1965, greenacre1984theory}.

\begin{property}\label{prop:transitionformula}
Transition formulas: 
The row factor scores can be computed from the 
column factor scores and vice versa.
\end{property}

CA/MCA/DiSCA/DiMCA analyze the rows and columns 
symmetrically; 
therefore, the row (respectively column) 
factors can be obtained from the data 
and the column (respectively row) 
factors by a transition formula.
Transition formulas can be derived 
from Equation \ref{gFactorScores_DiSCA}: 
by substituting the $ \DcNeg\gQ $ of $ \F $,
\begin{equation}\label{transition_knownF}
\F = \DrNeg\X\DcNeg\gQ, 
\quad\mathrm{with}\quad \DcNeg\gQ = \G\BDelta^{-1}
\end{equation}
and the $ \DrNeg\gP $ of $ \G $,
\begin{equation}\label{transition_knownG}
\G = \DcNeg\Xt \DrNeg\gP, \quad\mathrm{with}\quad \DrNeg\gP = \F\BDelta^{-1},
\end{equation}
$ \F $ can be computed from $ \G $ by
\begin{equation}\label{transition_F}
\F = \DrNeg\X\G\BDelta^{-1}
\end{equation}
and $ \G $ can be computed from $ \F $ by
\begin{equation}\label{transition_G}
\G = \DcNeg\Xt\F\BDelta^{-1}.
\end{equation}

\begin{property}\label{prop:supplementaryprojection}
Supplementary projections: 
Equations \ref{transition_F} and \ref{transition_G} 
can be used to estimate the factor scores 
from a supplementary, 
or called out-of-sample, 
row (respectively column) 
that is represented by the same set of columns (respectively rows).
\end{property}

The factor score of a supplementary row 
$\getc{\i}{\sup}$ (denoted $\getct{\fsup}{\sup}$) 
is computed with an equation similar to Equation \ref{transition_F}:
\begin{equation}
    \getct{\fsup}{\sup} = 
    \left(\getct{\i\transpose}{\sup}\1\right)^{-1}
    \getct{\i\transpose}{\sup}\G\BDelta^{-1}.
\end{equation}
The factor score of a supplementary column 
$\getc{\j}{\sup}$ (denoted $\getct{\gsup}{\sup}$ 
is computed by a similar equation to Equation \ref{transition_G}:
\begin{equation}
    \getct{\gsup}{\sup} = 
    \left(\1\transpose\getc{\j}{\sup}\right)^{-1}\getct{\j\transpose}{\sup}\F\BDelta^{-1}.
\end{equation}
Geometrically, these factor scores project 
the supplementary rows (or column) 
onto the component space built by the original data. 

\begin{property}\label{prop:distributionalequivalence}
Distributional equivalence: Identical rows (or columns) 
can be replaced by their sum without affecting the results.
\end{property}

Because CA and MCA analyze the frequencies of the occurrences, 
two rows (or two columns) that are proportional 
to each other become identical after the preprocessing steps. 
These identical rows (or columns) 
can be represented by \textit{two}
coincident points in the component space, 
and the two points can be merged into \textit{one} 
with the sum of the original weights
\citep{fichet2009metrics,benzecri1973analyse, greenacre1984theory}. 
In addition, merging these two points 
does not change the geometry of the component space.

\begin{property}\label{prop:barycentricprojection}
Barycentric projection: 
Row and column factor scores, 
 have barycenters of 0.
\end{property}

The row and the column factor scores of 
CA/MCA/DiSCA/DiMCA share a common barycenter 
(i.e., weighted mean) of 0; Formally
\begin{equation}\label{barycent_CA}
\frac{1}{I} \rr\transpose \getc{\f}{\ell} = \frac{1}{J} \cc\transpose \getc{\g}{\ell} = 0.
\end{equation}
Specifically in MCA, for each of the \textit{K} variables, 
the (column) factor scores of its levels will have a weighted mean of zero:
\jcy{\begin{equation}\label{barycent1_MCA}
\sum_{j = 1}^{J_k} \getsc{c}{j}{k} \g_{j, k, \ell} = 0,
\end{equation}
where $\getsc{c}{j}{k}$ is the column weight for the $j$th level of the $k$th variable, 
and $\g_{j, k, \ell}$ is the factor scores of the $j$th level of the $k$th variable 
on the $\ell$th component.}
In addition, for all observations that belong to this variable level, their mean (row) factor score will equal the (column) factor score of this variable level.

\begin{property}\label{prop:embeddedsolution}
The embedded solution: in CA, the GSVD of the non-centered matrix (i.e., $ \Z $) will have the first generalized singular value of 1, the first left generalized singular vector of $ \rr $, and the first right generalized singular vector of $ \cc $. In addition, the following components will be equivalent to the GSVD of the matrix $ \X $.
\end{property}

The embedded solution holds because the GSVD in CA (Equation \ref{CA_GSVD}) can be rewritten as:
\begin{equation}\label{rewriteGSVD_CA}
\begin{gathered}
\X = \Z-\rr\ct = \gP \BDelta \gQt = \sum_{\ell = 1}^L \getc{\delta}{\ell} \getc{\gp}{\ell} \getct{\gqt}{\ell} \\
\mathrm{under\;the\;constraints}\quad
\gPt\DrNeg\gP = \gQt\DcNeg\gQ = \eye
\end{gathered}
\end{equation}
which gives
\begin{equation}\label{noncentGSVD_CA}
\begin{gathered}
\Z = \rr\ct + \gP\BDelta \gQt =  1 \times \rr\ct +  \sum_{\ell = 2}^L \getc{\delta}{\ell} \getc{\gp}{\ell} \getct{\gqt}{\ell} \\
\mathrm{under\;the\;constraints}\quad
\gPt\DrNeg\gP = \gQt\DcNeg\gQ = \eye.
\end{gathered}
\end{equation}
Therefore, when the non-centered data $ \Z $ is analyzed, the first generalized singular value $\getc{\delta}{1}$ equals 1, the first left generalized singular vector $\getc{\gp}{1}$ equals $ \rr $, and the first right generalized singular vector $\getc{\gq}{1}$ equals $ \cc $. With $ \rr\ct $ computing the \textit{expected} frequencies of $ \Z $ under independence, the CA of $ \X $, where $ \X = \Z- \rr\ct$, analyzes the deviation of the \textit{observed} data (i.e., $ \Z $) from the independence (i.e., $ \rr\ct $).

\begin{property}\label{prop:asymmetricprojection}
The asymmetric projection: The row and column factor scores of each component can be scaled to have either a variance of 1 or a variance of the associated eigenvalue. When the row and the column factor scores are scaled differently (i.e., one to have a variance of 1 with the other having a variance of the eigenvalue), they are projected asymmetrically.
\end{property}

In the CA/MCA/DiSCA framework, the rows and columns are analyzed symmetrically and the extracted components can be seen from the perspective of the rows or of the columns. From the perspective of the rows, the component space is defined by the row factor scores (denoted by $ \widetilde{\F} $) which are computed as:
\begin{equation}\label{gAsymF_CA}
\widetilde{\F} = \DrNeg\gP
\quad\mathrm{with\;its\;variance\;}\quad
\widetilde{\F}\transpose\Dr\widetilde{\F} = \eye.
\end{equation}
The columns can then be projected onto this component space as $ \G $:
\begin{equation}\label{gSymG_CA}
\G = \DcNeg\gQ\BDelta
\quad\mathrm{with\;its\;variance\;}\quad
\Gt\Dc\G = \Blambda.
\end{equation}
Because the columns factor scores are projected onto this space and scaled differently from the row factor scores, this projection is \textit{asymmetric}.

From  the perspective of the columns, the component space is defined by the column factor scores (denoted by $ \widetilde{\G} $) which are computed as:
\begin{equation}\label{gAsymG_CA}
\widetilde{\G} = \DcNeg\gQ
\quad\mathrm{with\;its\;variance\;}\quad
\widetilde{\G}\transpose\Dc\widetilde{\G} = \eye.
\end{equation}
The asymmetric projection of the rows onto this space can then be computed as $ \F $:
\begin{equation}\label{gSymF_CA}
\F = \DrNeg\gP\BDelta
\quad\mathrm{with\;its\;variance\;}\quad
\Ft\Dr\F = \Blambda.
\end{equation}

According to Property \ref{prop:embeddedsolution}, because the first generalized singular value of the uncentered matrix equals 1, it is the maximum amount of variance any given component can have. When the factor scores of a component has a variance of 1, these factor scores span the entire space---called \textit{simplex}---of the data. From Property \ref{prop:barycentricprojection}, the simplex defined by the row and the simplex defined by the column factors share the same \textit{barycenter}. When one is used to define the simplex, the other can be projected asymmetrically onto this simplex where the distance between any two factor scores (including the distance between a row factor score and a column factor score) is meaningful. In contrast, when both the rows and the columns are symmetrically projected (as $ \F $ and $ \G $), only the distances within the same set are meaningful.



\newpage

\section{Transition formulas with and without sparsification}

We derive the transition formulas for row and column factor scores, for regular and sparse CA-related methods. We provide a step-by-step breakdown of these derivations, highlighting specifically where the projection operators introduce non-linearity, leading to sparsity.
 

The original transition formulas are derived as follows for row ($\f$) and column ($\g$) factor scores:
\begin{equation}
    \begin{split}
  \getc{\f}{\ell} = & \DrNeg\getc{\gp}{\ell}\getc{\delta}{\ell} \\
          = & \DrNeg(\gP\BDelta\gQt)\DcNeg\getc{\gq}{\ell} \\
          = & \DrNeg\X\DcNeg\getc{\gq}{\ell} \\
          = & \DrNeg(\Z-\rr\ct)\DcNeg\getc{\gq}{\ell} \\
          = & \DrNeg\Z\DcNeg\getc{\gq}{\ell} \\
          = & \R\DcNeg\getc{\gq}{\ell} \\
          = & \DrNeg\Z\getc{\g}{\ell}\getct{\delta^{-1}}{\ell},
    \end{split}
\end{equation}
where $\Z$ is the contingency table, $\X$ is the probability matrix, and $\R$ is the row profiles where $\R = \DrNeg\Z$ with each row of $\R$ sums to 1;
\begin{equation}
  \begin{split}
    \getc{\g}{\ell} = & \DcNeg\getc{\gq}{\ell}\getc{\delta}{\ell} \\
          = & \DcNeg(\gQ\BDelta\gPt)\DrNeg\getc{\gp}{\ell} \\
          = & \DcNeg\X\DrNeg\getc{\gq}{\ell} \\
          = & \DcNeg(\Z-\rr\ct)\transpose\DrNeg\getc{\gp}{\ell} \\
          = & \DcNeg\Zt\DrNeg\getc{\gp}{\ell} \\
          = & \C\DrNeg\getc{\gp}{\ell} \\
          = & \DcNeg\Zt\getc{\f}{\ell}\getct{\delta^{-1}}{\ell},
  \end{split}
\end{equation}
where $\C$ is the column profiles where $\C = \DcNeg\Zt$ with each column of $\C$ sums to 1.


With sparsification, the new transition formulas for CA/MCA/DiSCA require the projecting operators used in gsGSVD's algorithm: 
\begin{equation}
    \begin{split}
    \getc{\f}{\ell} &= \DrNeghalf\projLGLtwoOrth{\gp}{\P}{\DrNeghalf\X\getc{\g}{\ell}\getct{\delta^{-1}}{\ell}}\getc{\delta}{\ell} \\
    &= \DrNeghalf\projLGLtwoOrth{\gp}{\P}{\DrNeghalf(\Z - \rr\ct)\getc{\g}{\ell}\getct{\delta^{-1}}{\ell}}\getc{\delta}{\ell} \\
    &= \DrNeghalf\projLGLtwoOrth{\gp}{\P}{\DrNeghalf\Z\getc{\g}{\ell}\getct{\delta^{-1}}{\ell}}\getc{\delta}{\ell} \\
    &= \DrNeghalf\projLGLtwoOrth{\gp}{\P}{\DrNeghalf\Z\getc{\g}{\ell}\getct{\delta^{-1}}{\ell}}\getct{\delta}{\ell} \\
    &= \DrNeghalf\projLGLtwoOrth{\gp}{\P}{\Drhalf\R\getc{\gq}{\ell}}\getc{\delta}{\ell},
  \end{split}
\end{equation}

and
\begin{equation}
    \begin{split}
    \getc{\g}{\ell} &= \DcNeghalf\projLGLtwoOrth{\gq}{\Q}{\DcNeghalf\Xt\getc{\f}{\ell}\getct{\delta^{-1}}{\ell}}\getc{\delta}{\ell} \\
    &= \DcNeghalf\projLGLtwoOrth{\gq}{\Q}{\DcNeghalf(\Z -\rr\ct)\transpose\getc{\f}{\ell}\getct{\delta^{-1}}{\ell}}\getc{\delta}{\ell} \\
    &= \DcNeghalf\projLGLtwoOrth{\gq}{\Q}{\DcNeghalf\Zt\getc{\f}{\ell}\getct{\delta^{-1}}{\ell}}\getc{\delta}{\ell} \\
    &= \DcNeghalf\projLGLtwoOrth{\gq}{\Q}{\DcNeghalf\Zt\getc{\f}{\ell}\getct{\delta^{-1}}{\ell}}\getc{\delta}{\ell} \\
    &= \DcNeghalf\projLGLtwoOrth{\gq}{\Q}{\Dchalf\C\getc{\gp}{\ell}}\getc{\delta}{\ell}.
  \end{split}
\end{equation}

\newpage

\section{Simulated experiment to illustrate the loss of centering after projection}\label{append:centering}

Barycentric projection is based on the property that, when data are
(double) centered, any linear combination of the items will also be
centered. Which makes factor scores and loadings also centered, when
looking at the results of CA-related methods.

When sparsification is involved, however, this is no longer true,
because the way that the components and loadings are obtained is not
based on linear combinations anymore. However, there is one special case
where this property holds again: when applying
group-sparsification in the case of sparse MCA.

\begin{figure}[!!h]
    \centering
    \includegraphics[width=0.7\textwidth]{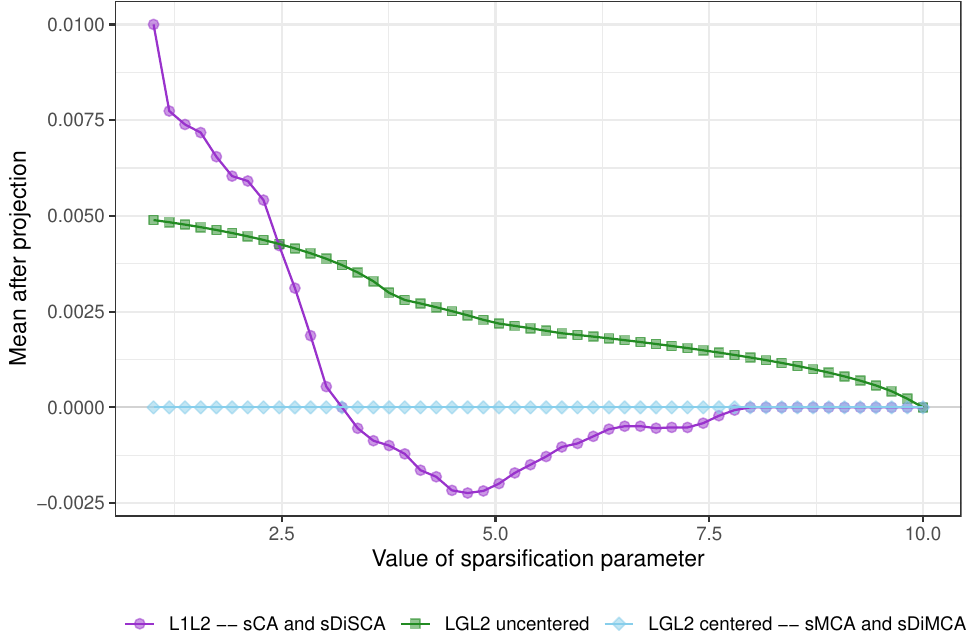}
    \caption{Projections onto the $\mathcal{L}_{1}$- or $\mathcal{L}_{\mathcal{G}}$- ball of a centered vector does not usually yield centered vectors, except in the case of the $\mathcal{L}_{\mathcal{G}}$ projection of a group-centered vector. Although the means of projections onto $\mathcal{L}_{1}$- and $\mathcal{L}_{2}$-balls (colored in purple) crossed 0, the mean is not exactly, but only close to, 0 and does not satisfy the barycentric property.}
    \label{fig:appcentering}
\end{figure}

To observe this property on a simple simulated example, let's consider a
normal random variable \(U\), and \(I = 100\) i.i.d. realisations of
this random variable, stored into a vector \(\mathbf{u}\). We then define the
vector \(\mathbf{x}\) of the centered values of \(\mathbf{u}\) and
\(\mathbf{x_\mathcal{G}}\) the vector of the ``group''-centered values
of \(\mathbf{u}\), where \(\mathcal{G}\) defines a partition of
\({1, \cdots, I}\) into 5 randomly assigned groups of 20.
We then apply the projection operators
\(\operatorname{proj}_{\mathcal{L}_1}\) and
\(\operatorname{proj}_{\mathcal{L_{\mathcal{G}}}}\) on these two vectors
with varying degrees of sparsity. The goal of this simulation is to
visualize the effect of these projections on centered data, especially
on the mean of the projected vectors, which is represented in Fig. \ref{fig:appcentering}.
We show on this plot the mean of the resulting
projected vectors, as a function of the value of the sparsity parameter.

First, we see that projecting a centered vector on an \(\mathcal{L}_1\)-
or an \(\mathcal{L}_{\mathcal{G}}\)-ball results in a vector that is
sparser, but not centered anymore. However, when the input vector is
group-centered, then the resulting (group-)projected vector is centered. This is
due to the fact that the \(\mathcal{L}_{\mathcal{G}}\)-projection
operator is scaling each group individually (therefore keeping each
group centered) or eliminating them.

\newpage

\section{Algorithms for plain SVD, GSVD, and sparse SVD (CSVD) and GSVD}\label{append:algorithms}

\setcounter{algocf}{0}
\renewcommand{\thealgocf}{E.\arabic{algocf}}

\begin{algorithm}
\LinesNotNumbered
 \DontPrintSemicolon
\KwData{$\X$, $\varepsilon$, $R$ \DontPrintSemicolon \Comment*[r]{Data $\X$, errors $\varepsilon$, and rank $R$}}
\KwResult{SVD of $\X$}
$\P \leftarrow \emptyset$; $\Q \leftarrow \emptyset$ \Comment*[r]{$\P$ and $\Q$ store the left and right singular vectors}
\For{$\ell=1, \dots, R$}{ 
  Initialize $\p^{(0)}$ and $\q^{(0)}$
  \Comment*[r]{Initialize $\p$ and $\q$}
  $\delta^{(0)} \leftarrow 0$;
  $\delta^{(1)} \leftarrow \p^{(0)}{\,}\transpose \X \q^{(0)}$
    \Comment*[r]{Update $\delta$ with $\p^{(0)}$ and $\q^{(0)}$}
  $t \leftarrow 0$\;
  \While{$\left(\Ltwo{\p^{(t+1)} - \p^{(t)}} \geq \varepsilon\right)$ \text{\rm or} $\left(\Ltwo{\q^{(t+1)} - \q^{(t)}} \geq \varepsilon\right)$}
  { $\p^{(t+1)} \leftarrow \X\q^{(t)}$\;
    $\q^{(t+1)} \leftarrow \X\transpose\p^{(t+1)}$\;
    $\delta^{(t+1)} \leftarrow \p^{(t+1)}{\,}\transpose \X \q^{(t+1)}$\;
    $t \leftarrow t+1$ \DontPrintSemicolon \Comment*[r]{Iterate until $\p$ and $\q$ are stable}
  }
  $\getc{\delta}{\ell} \leftarrow \delta^{(t)}$; $\getc{\p}{\ell} \leftarrow \p^{(t)}$;
  $\getc{\q}{\ell} \leftarrow \q^{(t)}$\;
  $\P \leftarrow \left[\,\P\,\mid\,\getc{\p}{\ell}\,\right]$\; 
  $\Q \leftarrow \left[\,\Q\,\mid\,\getc{\q}{\ell}\,\right]$\;
  $\X \leftarrow  \X - \delta_\ell\p_\ell\transpose\q_\ell$
  }
\caption{ALS algorithm of the SVD of $\X$}
\label{algo:svd}
\end{algorithm}

\begin{algorithm}
\LinesNotNumbered
\DontPrintSemicolon 
\KwData{$\X$, {\color{ForestGreen}$\M$, $\W$}, $\varepsilon$, $R$ \Comment*[r]{Data $\X$, {\color{ForestGreen} row metric matrix $\M$, column metric matrix $\W$}, error $\varepsilon$, and rank $R$}}
\KwResult{GSVD of $\X$}
$\widetilde{\X}$ = {\color{ForestGreen} $\Mhalf$}$\X${\color{ForestGreen} $\Whalf$}\; 
$\gP \leftarrow \emptyset$; $\gQ \leftarrow \emptyset$ \Comment*[r]{$\gP$ and $\gQ$ store the \textit{generalized} left and right singular vectors}
\For{$\ell=1, \dots, R$}{ 
  Initialize $\p^{(0)}$ and $\q^{(0)}$
  \Comment*[r]{Initialize $\p$ and $\q$}
  $\delta^{(0)} \leftarrow 0$;
  $\delta^{(1)} \leftarrow \p^{(0)}{\,}\transpose \widetilde{\X} \q^{(0)}$
    \Comment*[r]{Update $\delta$ with $\p^{(0)}$ and $\q^{(0)}$}
  $t \leftarrow 0$\;
  \While{$\left(\Ltwo{\p^{(t+1)} - \p^{(t)}} \geq \varepsilon\right)$ \text{\rm or} $\left(\Ltwo{\q^{(t+1)} - \q^{(t)}} \geq \varepsilon\right)$   
}  
  {
    $\p^{(t+1)} \leftarrow$ $\widetilde{\X}\q^{(t)}$\;
    $\q^{(t+1)} \leftarrow$ $\widetilde{\X}\transpose\p^{(t+1)}$\;
    $\delta^{(t+1)} \leftarrow \p^{(t+1)}{\,}\transpose \widetilde{\X} \q^{(t+1)}$\;
    $t \leftarrow t+1$ \Comment*[r]{Iterate until $\p$ and $\q$ are stable}
  }
  $\getc{\delta}{\ell} \leftarrow \delta^{(t)}$; $\getc{\p}{\ell} \leftarrow \p^{(t)}$; 
  $\getc{\q}{\ell} \leftarrow \q^{(t)}$\;
  {\color{ForestGreen} $\getc{\gp}{\ell}$} $\leftarrow$ {\color{ForestGreen} $\MNeghalf$}$\getc{\p}{\ell}$\;
  {\color{ForestGreen} $\getc{\gq}{\ell}$} $\leftarrow$ {\color{ForestGreen} $\WNeghalf$}$\getc{\q}{\ell}$\;
  $\gP \leftarrow \left[\,\gP\,\mid\,\getc{\gp}{\ell}\,\right]$\; 
  $\gQ \leftarrow \left[\,\gQ\,\mid\,\getc{\gq}{\ell}\,\right]$\; 
  $\widetilde{\X} \leftarrow  \widetilde{\X} - \delta_\ell\p_\ell\transpose\q_\ell$
  }
\caption{ALS algorithm of the GSVD of $\X$}
\algorithmfootnote{Note: The text colored in green is specific to the GSVD as compared to the SVD.}
\label{algo:gsvd}
\end{algorithm}

\begin{algorithm}
\LinesNotNumbered
\DontPrintSemicolon 
\KwData{$\X$, {\color{Maroon} $\getsc{\pntythes}{\p}{\ell}$}, {\color{Maroon} $\getsc{\pntythes}{\q}{\ell}$},  $\varepsilon$, $R$ \Comment*[r]{Data $\X$, errors $\varepsilon$, and rank $R$}}  \Comment*[r]{{\color{Maroon} sparse parameters $\getsc{\pntythes}{\p}{\ell}$ and $\getsc{\pntythes}{\q}{\ell}$ for singular vectors}} 
\KwResult{CSVD of $\X$ }
$\P \leftarrow \emptyset$; $\Q \leftarrow \emptyset$ \Comment*[r]{$\P$ and $\Q$ store the pseudo-left and right singular vectors}
\For{$\ell=1, \dots, R$}{ 
  Initialize $\p^{(0)}$ and $\q^{(0)}$
  \Comment*[r]{Initialize $\p$ and $\q$ either from SVD or randomly}
  $\delta^{(0)} \leftarrow 0$;
  $\delta^{(1)} \leftarrow \p^{(0)}{\,}\transpose \X \q^{(0)}$
    \Comment*[r]{Update $\delta$ with $\p^{(0)}$ and $\q^{(0)}$}
  $t \leftarrow 0$\;
  \While{$\left(\Ltwo{\p^{(t+1)} - \p^{(t)}} \geq \varepsilon\right)$ \text{\rm or} $\left(\Ltwo{\q^{(t+1)} - \q^{(t)}} \geq \varepsilon\right)$
  }{
    $\p^{(t+1)} \leftarrow$ {\color{Maroon} $\proj$(}$\X\q^{(t)}$, {\color{Maroon}$\Bs_{\Ler}(1) \cap \Bs_{\Lyi}(\getsc{\pntythes}{\p}{\ell}) \cap \P^\perp$)} \Comment*[r]{Projection of $\X\q^{(t)}$ onto the intersection}
    $\q^{(t+1)} \leftarrow$ {\color{Maroon} $\proj$(}$\X\transpose\p^{(t+1)}$, {\color{Maroon}$\Bs_{\Ler}(1) \cap \Bs_{\Lyi}(\getsc{\pntythes}{\q}{\ell}) \cap \Q^\perp$)}\;
    $\delta^{(t+1)} \leftarrow \p^{(t+1)}{\,}\transpose \X \q^{(t+1)}$\;
    $t \leftarrow t+1$ \Comment*[r]{Iterate until $\p$ and $\q$ are stable}
  }
  $\getc{\delta}{\ell} \leftarrow \delta^{(t)}$; $\getc{\p}{\ell} \leftarrow \p^{(t)}$; 
  $\getc{\q}{\ell} \leftarrow \q^{(t)}$\;
  $\P \leftarrow \left[\,\P\,\mid\,\getc{\p}{\ell}\,\right]$\; 
  $\Q \leftarrow \left[\,\Q\,\mid\,\getc{\q}{\ell}\,\right]$\;
}
\caption{ALS algorithm of the CSVD of $\X$}
\algorithmfootnote{Note: The text colored in red is specific to the CSVD as compared to the SVD.}
\label{algo:csvd}
\end{algorithm}

\begin{algorithm}
\LinesNotNumbered
\DontPrintSemicolon 
\KwData{$\X$, {\color{ForestGreen}$\M$, $\W$}, {\color{Maroon} $\getsc{\pntythes}{\p}{\ell}$}, {\color{Maroon} $\getsc{\pntythes}{\q}{\ell}$},  $\varepsilon$, $R$}  \Comment*[r]{Data $\X$, {\color{ForestGreen} row metric matrix $\M$, column metric matrix $\W$}, errors $\varepsilon$, and rank $R$} \Comment*[r]{{\color{Maroon} sparse parameters $\getsc{\pntythes}{\p}{\ell}$ and $\getsc{\pntythes}{\q}{\ell}$ for singular vectors}} 
\KwResult{sGSVD of $\X$ }
$\widetilde{\X}$ = {\color{ForestGreen} $\Mhalf$}$\X${\color{ForestGreen} $\Whalf$}\;
$\gP \leftarrow \emptyset$; $\gQ \leftarrow \emptyset$ \Comment*[r]{$\gP$ and $\gQ$ store the pseudo-generalized left and right singular vectors}
\For{$\ell=1, \dots, R$}{ 
  Initialize $\p^{(0)}$ and $\q^{(0)}$
  \Comment*[r]{Initialize $\p$ and $\q$ either from GSVD or randomly}
  $\delta^{(0)} \leftarrow 0$;
  $\delta^{(1)} \leftarrow \p^{(0)}{\,}\transpose \widetilde{\X} \q^{(0)}$
    \Comment*[r]{Update $\delta$ with $\p^{(0)}$ and $\q^{(0)}$}
  $t \leftarrow 0$
  \While{$\left(\Ltwo{\p^{(t+1)} - \p^{(t)}} \geq \varepsilon\right)$ \text{\rm or} $\left(\Ltwo{\q^{(t+1)} - \q^{(t)}} \geq \varepsilon\right)$
  }{
    $\p^{(t+1)} \leftarrow$ {\color{Maroon} $\proj$(} $\widetilde{\X}\q^{(t)}$,
    {\color{Maroon}$ \Bs_{\Lyi}(\getsc{\pntythes}{\p}{\ell}) \cap \Bs_{\Ler}(1) \cap \P^\perp$)}\;
    $\q^{(t+1)} \leftarrow$ {\color{Maroon} $\proj$(} $\widetilde{\X}\transpose\p^{(t+1)}$,
    {\color{Maroon}$ \Bs_{\Lyi}(\getsc{\pntythes}{\q}{\ell}) \cap \Bs_{\Ler}(1) \cap \Q^\perp$)}\;
    $\delta^{(t+1)} \leftarrow \p^{(t+1)}{\,}\transpose \widetilde{\X} \q^{(t+1)}$\;
    $t \leftarrow t+1$ \Comment*[r]{Iterate until $\p$ and $\q$ are stable}
  } 
  $\getc{\delta}{\ell} \leftarrow \delta^{(t)}$; $\getc{\p}{\ell} \leftarrow \p^{(t)}$; 
  $\getc{\q}{\ell} \leftarrow \q^{(t)}$\;
  {\color{ForestGreen} $\getc{\gp}{\ell}$} $\leftarrow$ {\color{ForestGreen} $\MNeghalf$}$\getc{\p}{\ell}$\;
  {\color{ForestGreen} $\getc{\gq}{\ell}$} $\leftarrow$ {\color{ForestGreen} $\WNeghalf$}$\getc{\q}{\ell}$\;
  $\gP \leftarrow \left[\,\gP\,\mid\,\getc{\gp}{\ell}\,\right]$\; 
  $\gQ \leftarrow \left[\,\gQ\,\mid\,\getc{\gq}{\ell}\,\right]$\;
}
\caption{ALS algorithm of sGSVD of $\X$}
\algorithmfootnote{Note: The text colored in red 
describes the sparsification constraints 
of the CSVD that are also used in the sGSVD, 
and the text colored in green describes 
the metric constraints of the 
GSVD that are also used in the sGSVD.}
\label{algo:sGSVD}
\end{algorithm}

\newpage
\section{Supplementary Figures}

\setcounter{figure}{0}

\begin{figure}
  \begin{subfigure}{\linewidth}
    \begin{center}
    \includegraphics[width=0.75\textwidth]{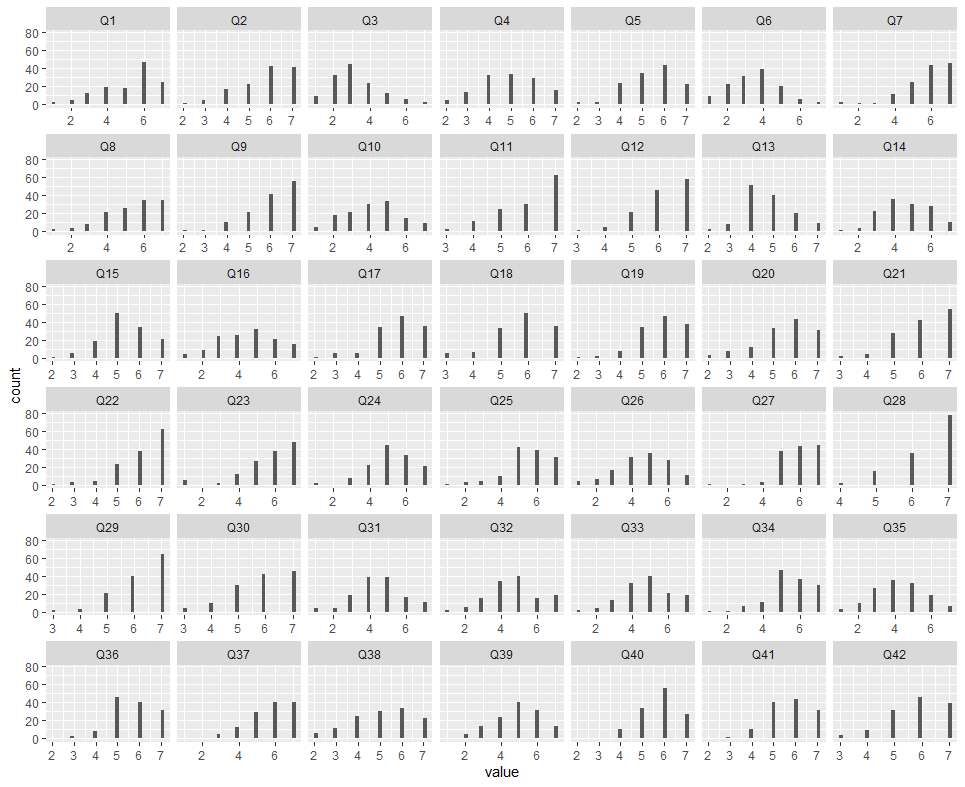}
    \end{center}
    \caption{The histogram of each item.}
    \label{sMCABin:a}
  \end{subfigure}
  \begin{subfigure}{\linewidth}
    \begin{center}
    \includegraphics[width=0.75\textwidth]{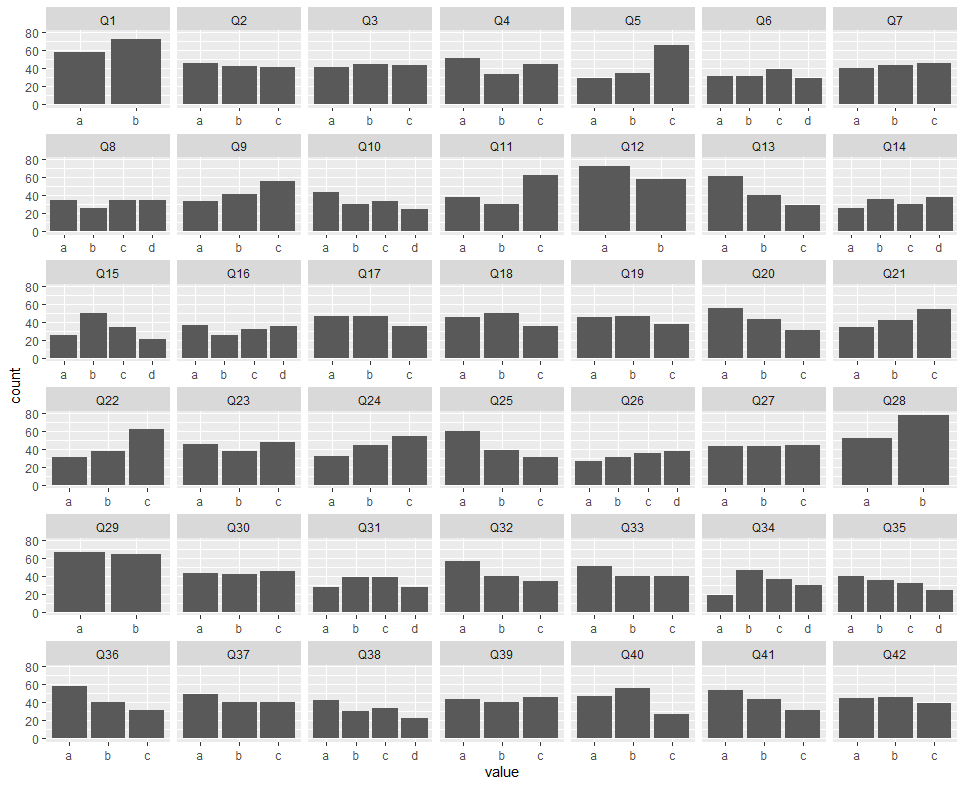}
    \end{center}
    \caption{The histogram of each item after binning.}
    \label{sMCABin:b}
  \end{subfigure}
  \caption{The grouping of item responses.}
  \label{sMCABin}
\end{figure}

\begin{figure}
    \caption{Pattern of the associations between items. The items are grouped by hierarchical clustering analysis with minimized Ward distances.}
    \label{sMCAheat}
    \centering
    \includegraphics[width=\textwidth]{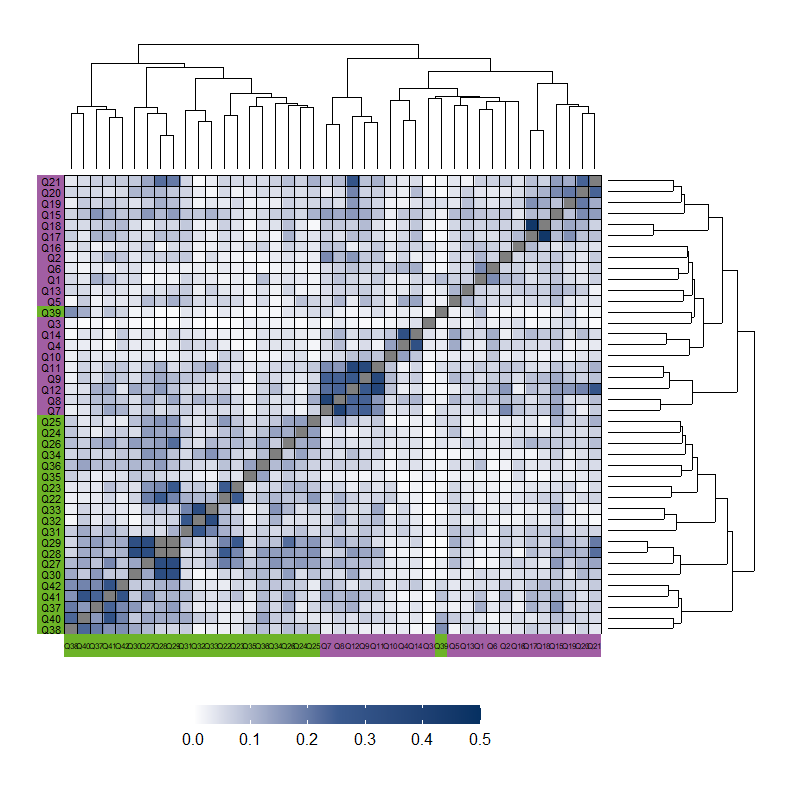}
\end{figure}

\end{document}